\documentclass[final]{siamart251216}
\usepackage[utf8]{inputenc}

\newsiamremark{remark}{Remark}

\usepackage{amssymb}
\usepackage{amsmath}
\usepackage{algorithm}
\usepackage{algorithmic}
\usepackage{color}
\usepackage{soul}
\usepackage{adjustbox}
\usepackage{placeins}
\usepackage{graphicx}
\usepackage{subcaption}
\usepackage{booktabs}

\newcommand\Tstrut{\rule{0pt}{2.6ex}}       
\newcommand\Hstrut{\rule{0pt}{3.0ex}} 

\title{Backward Error and Scaling in Level-Based \\ILU($\lowercase{k}$) and IC($\lowercase{k}$) Preconditioners}
\author{Jennifer Scott\thanks{Scientific Computing Department,
STFC Rutherford Appleton Laboratory, Didcot, Oxfordshire OX11 0QX, UK and
Department of Mathematics and Statistics, University of Reading, Whiteknights,
Reading RG6 6AX, UK (\email{jennifer.scott@stfc.ac.uk})}
\and
Miroslav T\r{u}ma\thanks{Charles University, Prague, Czech Republic
  (\email{mirektuma@karlin.mff.cuni.cz})}
}

\begin{document}

\maketitle

\begin{abstract}
We analyse the effect of rounding errors and diagonal scaling on level-based
ILU($k$) and IC($k$) preconditioners. Because the sparsity pattern is fixed in
advance, the computed factors are those that the incomplete factorization
produces, in exact arithmetic, for a nearby matrix, and the deviation of the
preconditioned matrix from the identity is bounded by an incompleteness term
plus a rounding term. When no shift or pivot modification is needed, scaling
leaves the iterates unchanged in exact arithmetic but can change the quantities
in this bound by orders of magnitude. For IC($k$), we bound how far rounding moves
the preconditioned eigenvalues from those of the incomplete factorization of
the nearby matrix; the bound involves Skeel condition numbers of the incomplete
factor and does not depend on the scaling. We also show that after Jacobi
scaling the entries of the factor are
bounded, so that in half precision a breakdown can be detected before any
overflow. Experiments in double, single and half precision illustrate the theory.
\end{abstract}

\begin{keywords}
level-based incomplete factorization, incomplete Cholesky, preconditioning,
backward error analysis, Skeel condition number, mixed precision
\end{keywords}

\begin{MSCcodes}
65F08, 65F10, 65F50, 65F35, 65G50
\end{MSCcodes}

\section{Introduction}

Our interest is in solving large sparse linear systems $Ax = b$
using preconditioned Krylov subspace methods.
Level-based incomplete factorization preconditioners ILU($k$) (and IC($k$) for symmetric
positive definite (SPD) systems), introduced by Watts in 1981~\cite{watt:81}, are
particularly well suited to sparse systems arising from the
discretisation of partial differential equations, where the
sparsity pattern is structured and the matrix is often
diagonally dominant or symmetric positive definite.
ILU($k$) remains widely used for these problems because it is simple, its fill
is adjustable, it can be computed efficiently in parallel, and it suits sparse
data structures \cite{hypo:01,meva:77,saad:03}.
However, in common with other incomplete factorization approaches,
level-based factorizations are sensitive to the ordering of $A$ and to the numerical behaviour of the elimination process~\cite{dume:89,saad:03}.
The fill level $k$ trades the cost of computing and applying the
preconditioner against its quality. It is usually small, since the fill can
grow rapidly with $k$~\cite{benz:02}, but larger values can be useful in
parallel~\cite{hypo:01}.
The quality of such preconditioners can vary substantially for
different problems and different $k$.

We develop a backward error analysis of ILU($k$) and IC($k$). It uses
only the fact that the sparsity pattern of the factors is fixed by the structure
of $A$ before any numerical computation, so the results hold for any incomplete
factorization with a pattern chosen this way, but not for threshold-based
methods. The analysis covers factors computed in low precision and applied
within Krylov iterations in high precision~\cite{sctu:2024}, and adaptive
precision factorizations, which compute in one precision and store different
entries of the factors in different precisions~\cite{lims:2026}.

The computed factors are those that the incomplete factorization produces, in
exact arithmetic, for a matrix that differs from $A$ by a backward error on the
sparsity pattern. The deviation of the preconditioned matrix from the identity
is then at most an incompleteness term plus a rounding term, the backward error times the norm of the inverse of the
preconditioner. Both terms are formed from the computed factors, so rounding
errors also affect the first.

In our experiments, $A$ is first scaled to $\widehat A$. In exact arithmetic this
leaves the incomplete factorization unchanged: the preconditioner and the
iterates are the same as without scaling. This no longer holds when a small
pivot, which can occur even when $A$ is well conditioned, is replaced by a
larger value, as in our ILU($k$) code, or when IC($k$) breaks down and is
repeated for the shifted matrix $\widehat A + \alpha I$ with $\alpha > 0$. When
no pivot is replaced and no shift is used, the norms in the bounds (of
the factors, of the inverse of the preconditioner and of the deviation from the
identity) depend on the scaling: with a different choice, they can change by
orders of magnitude. These norms therefore reflect the scaling as well as the
preconditioner, and in our experiments none of them predicts the iteration
count of the Krylov solver.

For IC($k$), bounding the backward error by the size of the factor gives a
bound on the rounding term that is at least the unit roundoff times the
condition number of the preconditioner, so it is uninformative when the
preconditioner is ill conditioned relative to the working precision. We
therefore also bound how far rounding errors move the preconditioned
eigenvalues, which determine CG convergence.

The main results are as follows. In Section~\ref{sec:main}, we show that the
backward error of ILU($k$) is at most the product of the absolute values of the
computed factors times a constant that depends on the largest number of entries
in a row or column of the factors, not on $n$: because the pattern is fixed in
advance, each entry of the factors is computed by a single formula with at most
that many terms. The deviation is then bounded by
the incompleteness and rounding terms. In Section~\ref{sec:ic}, for IC($k$),
we bound how far rounding moves the preconditioned eigenvalues from those that
IC($k$) gives, without rounding, for the nearby matrix. The bound is a small
multiple of the unit roundoff times the product of two Skeel condition numbers
of the incomplete factor; unlike the normwise bound, it has no lower limit
proportional to the condition number of the preconditioner, and when no shift
is used it does not depend on the scaling. We also show that, after a diagonal scaling that
makes the diagonal entries at most $1$, such as Jacobi scaling, the entries of
the factor are bounded, so that in half precision a breakdown is
detected before any overflow can occur. Numerical results for ILU($k$) and IC($k$) are reported in
Sections~\ref{sec:LU examples} and~\ref{sec:experiments}.
Section~\ref{sec:background} reviews related work and gives the notation and
the scaling convention, and Section~\ref{sec:conclusions} concludes.
\section{Background and Preliminaries}
\label{sec:background}

Given a sparse matrix $A \in \mathbb{R}^{n \times n}$ and a fill
level $k \geq 0$, the ILU($k$) factorization performs Gaussian elimination and
discards any fill entry with a \emph{level of fill} greater than $k$. The positions of
the entries of $A$, and the diagonal, have level $0$, and all other positions
have level $\infty$. When eliminating with pivot $s$ updates position $(i,j)$
from the entries in positions $(i,s)$ and $(s,j)$, the level of $(i,j)$ becomes
\[
  \mathrm{lev}(i,j) \leftarrow
  \min\{\mathrm{lev}(i,j),\ \mathrm{lev}(i,s) + \mathrm{lev}(s,j) + 1\},
\]
the usual rule, which our implementations use. ILU($k$) keeps the positions
with level at most $k$; ILU(0) keeps only the positions of the entries of $A$.

\subsection{Related work}
\label{sec:related}

Classical work on ILU methods focuses on breakdown avoidance,
ordering effects, fill control, and the influence of sparsity
structure on robustness; see, for example,
\cite{dume:89,meva:77} and the many subsequent developments
of modified ILU, pivoting and shifting strategies \cite{limo:99,mant:80}.
Kershaw~\cite{kers:78} gives an early example
of unexpected breakdown in incomplete Cholesky, and Chow and Saad~\cite{chsa:97}
present an experimental
study showing that ILU preconditioners can fail even without
small pivots, and that failure is not reliably predicted by any
single diagnostic. They identify inaccuracy and instability of the triangular
solves as major causes of failure, and discuss how scaling, reordering and
diagonal perturbation can sometimes mitigate them.

The finite-precision behaviour of incomplete factorizations has
been analysed, generally assuming a single working precision throughout
the factorization and focusing either on ILU(0) or on factorizations
with numerical dropping.
Bollh\"ofer~\cite{boll:03} shows that for threshold-based incomplete factorizations the quality of the
preconditioned operator can be monitored and partially controlled by the norms of the inverse
triangular factors, and develops a dropping strategy that keeps
these bounded during the factorization. Bru, Mar\'{\i}n, Mas and
T\r{u}ma~\cite{bmmt:08} propose balancing the growth in ILU and in the
computed inverse triangular factors.
The work of these authors applies to value-dependent factorizations in which the
sparsity pattern evolves during elimination; the behaviour of  ILU($k$), where the sparsity pattern
is fixed in advance of the numerical factorization, has received
less attention, although there are important theoretical and practical contributions by Hysom and Pothen~\cite{hypo:01, hypo:02}.

Low-precision arithmetic is faster and more energy efficient than double
precision on current hardware~\cite{hima:2022}.
Much of the recent theory concerns iterative refinement based on the
generalized minimal residual method (GMRES-IR), where low-precision factors are used within
an inner Krylov solve while the outer refinement iteration is
performed in high precision~\cite{cahi:2018,hipr:2021}. Amestoy et
al.~\cite{ames:2023} model an approximate sparse factorization as the exact
factorization of a perturbed matrix and separate the approximation and rounding
errors in the analysis of GMRES-IR, and the stability of split-preconditioned
FGMRES in several precisions is analysed in~\cite{cada:2024}. A survey of
mixed-precision methods is given in~\cite{hima:2022}.
Scott and T\r{u}ma~\cite{sctu:2024} show how breakdown of incomplete
factorizations in half precision can be avoided by prescaling $A$, restarting
with a larger global shift after a breakdown, and using safe tests to detect
possible overflow. They illustrate this with level-based IC($k$) within
iterative refinement based on CG or GMRES. In their experiments, restarting after a breakdown always
avoided overflow, and the iteration count with the half precision factor may decrease little as
$k$ increases. Lemma~\ref{lem:bounded-entries} explains the first observation, and the
results of Section~\ref{sec:experiments}, with Jacobi scaling, suggest the
reason for the second.

We consider ILU($k$) and
IC($k$) applied within a Krylov method.
The behaviour of the preconditioner depends on the growth and
conditioning of the incomplete factors as well as on the residual: incomplete
factorizations may produce factors with significant growth or small pivots,
giving ill-conditioned triangular systems in which both the residual and the
finite-precision perturbations are amplified \cite{boll:03,chsa:97,elma:86,saad:03}.

\subsection{Notation and floating-point model}
\label{sec:notation}
For $B \in \mathbb{R}^{n \times n}$, the spectral and Frobenius norms
are denoted by $\|B\|_2$ and $\|B\|_F$, respectively.
If no subscript is given, $\|\cdot\| := \|\cdot\|_2$. For nonsingular $B$,
$\kappa(B) := \|B\|_2\,\|B^{-1}\|_2$, also written $\kappa_2(B)$, and
$\kappa_\infty(B) := \|B\|_\infty\|B^{-1}\|_\infty$.

The incomplete factorization (ILU($k$) or IC($k$)) is computed and stored in a
working precision with
unit roundoff $u$; in our experiments this is IEEE half (fp16), single (fp32)
or double (fp64) precision, with
$u = 2^{-11}$, $2^{-24}$ and $2^{-53}$ respectively. Matrix
scaling, the Krylov iteration and all applications of $A$, $A^T$ and the
preconditioner are performed in fp64, which we treat as exact. For the
factorization we use the standard floating-point model: for floating-point
numbers $a$ and $b$,
\[
  \mathrm{fl}(a \;\mathrm{op}\; b)
  =
  (a \;\mathrm{op}\; b)(1+\delta),
  \qquad
  |\delta| \le u,
  \qquad
  \mathrm{op}\in\{+,-,\times,/\},
\]
together with the standard quantity
\[
  \gamma_p = \frac{p\,u}{1-p\,u},
\]
which describes the accumulation of rounding errors over $p$
floating-point operations. The model ignores underflow. All codes are compiled
with the NAG Fortran compiler, release~7.1, with the options
\texttt{-round\_hreal}, so that every operation on fp16 numbers is rounded to
fp16, and \texttt{-ieee=full}, so that underflow is gradual (subnormal numbers
are kept). The BLAS routines are the original Fortran
implementations~\cite{blas2,blas1}.

We write $L_k$ and $U_k$ for the unit lower triangular and upper triangular
ILU($k$) factors of $A$ in exact arithmetic, $\widetilde L_k$ and
$\widetilde U_k$ for the factors computed in floating-point arithmetic, and
\begin{equation}
  M_k := L_kU_k,
  \qquad
  \widetilde{M}_k := \widetilde{L}_k \widetilde{U}_k
  \label{eq:M-def}
\end{equation}
for the corresponding preconditioners. The sparsity patterns $\mathcal{S}_L$
and $\mathcal{S}_U$ of the ILU($k$) factors are determined by the structure of
$A$ and $k$ alone~\cite{meva:77,saad:03}, so the computed factors have the
same patterns as the factors in exact arithmetic.

In exact arithmetic, the incomplete factorization satisfies
\begin{equation}
  A = L_kU_k + R_k,
  \label{eq:ILU-decomp}
\end{equation}
 where $R_k$ is the
incomplete factorization residual consisting of the discarded fill.

If $\widetilde M_k$ is used in stationary iterative refinement, with the
residual computed in higher precision, the iteration converges if
$\|A\widetilde M_k^{-1} - I\|_2 < 1$~\cite[Thm.~3.1]{cahi:2018}. For GMRES with
right preconditioning, choosing the polynomial $(1-z)^j$ in the minimization
gives $\|r_j\|_2 \le \|A\widetilde M_k^{-1} - I\|_2^{\,j}\,\|r_0\|_2$, where
$r_j = b - Ax_j$, so
$\|A\widetilde M_k^{-1} - I\|_2 < 1$ is also sufficient for GMRES to converge,
but it is not necessary.

\subsection{Scaling}
\label{sec:scaling}
In all our numerical experiments $A$ is scaled before the incomplete
factorization is computed. We write $\widehat{A} = D_LAD_R$ for the scaled matrix
and use a hat for quantities associated with it; in the symmetric positive
definite case the scaling is symmetric, $\widehat{A} = DAD$ with
$D = \mathrm{diag}(A)^{-1/2}$ (Jacobi scaling), so that $\widehat{A}$ has unit
diagonal.
For the nonsymmetric problems, $D_L$ and $D_R$ are the row and column scaling
factors of {\tt MC77}~\cite{Ruiz:2001}, which equilibrates the rows and columns
in the $\infty$-norm. The factorization is applied to $\widehat{A}$, giving $\widehat{M}_k = \widehat{L}_k\widehat{U}_k$,
but the Krylov method is applied to the original system $Ax = b$, so that the
stopping criterion refers to the original residual. The scaling is undone at
each application of the preconditioner, so the operator applied is
$D_R\,\widetilde M_k^{-1}D_L$, with $\widetilde M_k$ computed from $\widehat A$.

In exact arithmetic, and provided no shift or pivot modification is used, the
operator $D_R\,\widehat{M}_k^{-1}D_L$ is $M_k^{-1}$, the ILU($k$) preconditioner
of $A$ itself, since the factors of
$\widehat A$ are $\widehat L_k = D_LL_kD_L^{-1}$ and $\widehat U_k = D_LU_kD_R$
($\widehat L_k = DL_k$ for IC($k$)); prescaling therefore leaves the
preconditioner, and hence the iterates, unchanged, as is well known (see,
e.g.,~\cite{limo:99}). Three things do depend on the scaling. The factors
are computed in finite precision, so the rounding errors they incur are those of
$\widehat{A}$ and not of $A$. If a shift is used
(Section~\ref{sec:ic-positivity}), the factorization is of
\[
  \widehat A + \alpha I = D\,(A + \alpha D^{-2})\,D,
\]
so the shift is relative to the diagonal of $A$ and depends on $D$: with Jacobi
scaling, $a_{ii}$ is replaced by $(1+\alpha)\,a_{ii}$. Pivot modification
(Section~\ref{sec:backward-error}) also depends on the scaling, since the
replacement is formed from the entries of $\widehat A$.

The results of Sections~\ref{sec:main} and~\ref{sec:ic} are stated for a general
matrix and are applied with $\widehat{A}$ in place of $A$: every norm, growth factor
and condition number appearing in the bounds is that of $\widehat{A}$ and its
factors. Since $AM_k^{-1} = D_L^{-1}(\widehat A\widehat M_k^{-1})D_L$ and
$M_k^{-1}A = D_R(\widehat{M}_k^{-1}\widehat{A})D_R^{-1}$, the preconditioned
matrices for $A$ and $\widehat A$ are similar: they have the same spectrum, but
their norms can differ by as much as $\kappa(D_L)$ or $\kappa(D_R)$.

The analysis needs only that the patterns depend on the structure of $A$
alone. ILU($k$) and IC($k$) are such factorizations, with nested patterns
indexed by $k$~\cite{hypo:02}. Threshold-based methods such as ILUT and
ILU($\tau$) are not: their pattern depends on intermediate values, so rounding
errors and scaling can change it.

\section{Backward error analysis of ILU($k$)}
\label{sec:main}

\subsection{Backward error of the computed ILU(\textit{k}) factors}
\label{sec:backward-error}

\begin{lemma}[Backward error of the computed ILU($k$) factors]
\label{lem:backward-error}
Let the ILU($k$) factorization of $A$ be computed in floating-point arithmetic
with unit roundoff $u$, with no zero pivot, giving $\widetilde L_k$ and
$\widetilde U_k$. Let $\mathcal{S} := \mathcal{S}_L \cup \mathcal{S}_U$, which
contains the pattern of $A$, let $m_k$ be the maximum number of positions in a
row of $\mathcal S_L$ or a column of $\mathcal S_U$ (including the diagonal), and assume $m_ku < 1$.
For a matrix $B$,
let $B|_{\mathcal S}$ have the entries $b_{ij}$ with $(i,j) \in \mathcal S$ and
zeros elsewhere, and let $B|_{\mathcal S^c} := B - B|_{\mathcal S}$. Define
\begin{equation}
  \widetilde R_k := A - \widetilde L_k\widetilde U_k, \qquad
  \Delta A_k := \widetilde R_k|_{\mathcal S}, \qquad
  E_k := \widetilde R_k|_{\mathcal S^c}.
  \label{eq:DA-def}
\end{equation}
Then
\begin{enumerate}
\item[\textup{(i)}]
  $|\Delta A_k| \le \gamma_{m_k}\,
  \bigl(|\widetilde L_k|\,|\widetilde U_k|\bigr)\big|_{\mathcal S}$,
  where $|\cdot|$ and the inequality are understood entrywise;
\item[\textup{(ii)}]
  $\widetilde L_k$ and $\widetilde U_k$ are the factors that the incomplete
  factorization with pattern $\mathcal S$ produces, in exact arithmetic, for
  $\widetilde A := A - \Delta A_k$, and $E_k = \widetilde A -
  \widetilde L_k\widetilde U_k$ is its residual.
\end{enumerate}
\end{lemma}

\begin{proof}
The pattern is fixed symbolically before any numerical computation, so ILU($k$)
is Gaussian elimination in which every update to a position outside
$\mathcal S$ is discarded. For $(i,j) \in \mathcal S$ the computed entry is
therefore
\[
  \tilde u_{ij} = \mathrm{fl}\Bigl(a_{ij} - \sum_{s<i}\tilde\ell_{is}\tilde u_{sj}\Bigr)
  \ (i \le j),
  \qquad
  \tilde\ell_{ij} = \mathrm{fl}\Bigl(\bigl(a_{ij} - \sum_{s<j}\tilde\ell_{is}\tilde u_{sj}\bigr)
  \big/\tilde u_{jj}\Bigr)
  \ (i > j),
\]
where the sums are over $s$ with $(i,s), (s,j) \in \mathcal S$ and are
evaluated in some order. The standard result on such
expressions~\cite[Lem.~8.4]{high:2002} states that if
$y = \bigl(c - \sum_{s=1}^{p-1} a_sb_s\bigr)/b_p$ is evaluated in
floating-point arithmetic, in any order, then the computed $\tilde y$ satisfies
\[
  \Bigl|c - \sum_{s=1}^{p-1} a_sb_s - b_p\tilde y\Bigr|
  \le \gamma_p\Bigl(\sum_{s=1}^{p-1}|a_sb_s| + |b_p\tilde y|\Bigr),
\]
with $\gamma_{p-1}$ in place of $\gamma_p$ if $b_p = 1$, so that there is no
division. We apply it with $c = a_{ij}$, $a_sb_s = \tilde\ell_{is}\tilde u_{sj}$, and
$b_p\tilde y = \tilde\ell_{ij}\tilde u_{jj}$ if $i > j$ or
$b_p\tilde y = \tilde\ell_{ii}\tilde u_{ij} = \tilde u_{ij}$ if $i \le j$. The
number of terms $p$ is at most the number of positions in row $i$ of
$\mathcal S_L$ and in column $j$ of $\mathcal S_U$, so $p \le m_k$, and
\[
  \Bigl|a_{ij} - \sum_{s \le \min(i,j)} \tilde\ell_{is}\tilde u_{sj}\Bigr|
  \le \gamma_{m_k} \sum_{s \le \min(i,j)} |\tilde\ell_{is}|\,|\tilde u_{sj}| ,
\]
which is (i). For (ii), $\widetilde A - \widetilde L_k\widetilde U_k
= \widetilde R_k - \Delta A_k = E_k$, which is zero in every position of
$\mathcal S$. With nonzero pivots, the factors with pattern $\mathcal S$ are
determined uniquely by the condition that their product agrees with
$\widetilde A$ in every position of $\mathcal S$, which proves (ii).
\end{proof}

In (ii) the pattern is $\mathcal S$, the ILU($k$) pattern of $A$, not that of
$\widetilde A$. The two can differ: $\Delta A_k$ is generally nonzero at fill
positions, where $A$ has zeros, and the symbolic ILU($k$) rules applied to
$\widetilde A$ would treat these positions as entries of level $0$ and could
admit further fill. The matrix $\widetilde A$ appears only in the analysis; it
is never formed. In the fp32 and fp16 runs of
Section~\ref{sec:experiments}, $\widehat A$ is rounded to the working precision
before the factorization. This adds at most $u|\mathrm{fl}(\widehat A)|$ to $\Delta A_k$,
and since $\gamma_{m_k} + u(1+\gamma_{m_k}) \le \gamma_{m_k+1}$,
Lemma~\ref{lem:backward-error} then holds for $\widehat A$ with
$\gamma_{m_k+1}$ in place of $\gamma_{m_k}$.
The constant in Lemma~\ref{lem:backward-error} is $\gamma_{m_k}$, compared
with $\gamma_n$ for Gaussian elimination~\cite[Thm.~9.3]{high:2002}.
In an adaptive precision factorization~\cite{lims:2026}, the operations are
performed with unit roundoff $u$ and each entry is stored in a format with unit roundoff
$u_{ij}$; later operations use the stored values. Each equation in the proof
defines one new entry, so its storage adds one rounding, and (i) holds with
$\gamma_{m_k}$ replaced at position $(i,j)$ by approximately
$\gamma_{m_k} + u_{ij}$.

ILU($k$) does not incorporate pivoting, although $A$ may be preordered; the ordering is omitted
from the notation. A small pivot gives inaccurate factors, particularly in low
precision.
Some ILU implementations address this by replacing any pivot of magnitude below a
threshold $\tau$ by a larger value. Let $D_q$ be the diagonal matrix of the changes
made to the pivots. Since each replacement increases the magnitude of the
pivot, Lemma~\ref{lem:backward-error} holds with $A + D_q$ in place of $A$,
and Theorem~\ref{thm:stability} then bounds $\|A\widetilde M_k^{-1} - I\|_2$
with the extra term $\|D_q\widetilde M_k^{-1}\|_2$.
In our code, following Li and Shao~\cite{lish:2011}, a pivot with
$|\tilde u_{qq}| < \tau = 10^{-10}$ is replaced by
$\beta\max_i|\hat a_{iq}|$, where $\hat a_{ij}$ are the entries of $\widehat A$, with $\beta = 10^{-2(1-q/n)}$, which increases with
the column index $q$.

\subsection{Deviation of the computed preconditioner}
\label{sec:main-theorem}
In exact arithmetic, $A = M_k + R_k$ by~\eqref{eq:ILU-decomp}, and the
deviation of the preconditioned operator from the identity is
\begin{equation}
  \eta_k := \|AM_k^{-1} - I\|_2 = \|R_kM_k^{-1}\|_2 .
  \label{eq:eta-def}
\end{equation}

\begin{theorem}[Deviation bound for ILU($k$)]
\label{thm:stability}
Suppose that the assumptions of Lemma~\ref{lem:backward-error} hold, and let
$\widetilde M_k = \widetilde L_k\widetilde U_k$, with $\Delta A_k$ and $E_k$
as in~\eqref{eq:DA-def}. Then
\[
  A\widetilde M_k^{-1} - I = E_k\widetilde M_k^{-1} + \Delta A_k\widetilde
  M_k^{-1},
\]
and hence
\begin{equation}
  \|A\widetilde{M}_k^{-1} - I\|_2
  \le \tilde\eta_k + \omega_k,
  \qquad
  \tilde\eta_k := \|E_k\widetilde M_k^{-1}\|_2,
  \quad
  \omega_k := \|\Delta A_k\widetilde M_k^{-1}\|_2 ,
  \label{eq:main-bound-revised}
\end{equation}
with
\begin{equation}
  \omega_k
  \le \gamma_{m_k}\,\bigl\||\widetilde L_k|\,|\widetilde U_k|\bigr\|_2\,
      \|\widetilde M_k^{-1}\|_2
  \le \gamma_{m_k}\,\|\widetilde L_k\|_F\,\|\widetilde U_k\|_F\,
      \|\widetilde M_k^{-1}\|_2 .
  \label{eq:rho-bound}
\end{equation}
\end{theorem}

\begin{proof}
$\widetilde M_k$ is nonsingular because $\widetilde L_k$ is unit lower
triangular and the pivots are nonzero. By~\eqref{eq:DA-def},
$A - \widetilde M_k = \widetilde R_k = E_k + \Delta A_k$; multiplying on the
right by $\widetilde M_k^{-1}$ gives the identity, and taking norms
gives~\eqref{eq:main-bound-revised}. For~\eqref{eq:rho-bound}, recall that
$\|B\|_2 \le \||B|\|_2$ for any $B$, and $\|B\|_2 \le \|C\|_2$ if
$0 \le B \le C$ entrywise. By Lemma~\ref{lem:backward-error}(i),
$|\Delta A_k| \le \gamma_{m_k}|\widetilde L_k||\widetilde U_k|$ entrywise, so
\[
  \omega_k \le \|\Delta A_k\|_2\|\widetilde M_k^{-1}\|_2
  \le \gamma_{m_k}\bigl\||\widetilde L_k||\widetilde U_k|\bigr\|_2
      \|\widetilde M_k^{-1}\|_2 .
\]
The second inequality in~\eqref{eq:rho-bound} follows from
\[
  \bigl\||\widetilde L_k||\widetilde U_k|\bigr\|_2 \le
  \bigl\||\widetilde L_k||\widetilde U_k|\bigr\|_F \le
  \|\widetilde L_k\|_F\|\widetilde U_k\|_F .
\]
\end{proof}

The bound~\eqref{eq:main-bound-revised}  splits the residual
$\widetilde R_k$ of the computed factors into $E_k$, outside the pattern, which
is due to the incompleteness of the factorization, and $\Delta A_k$, on the
pattern, which is due to rounding errors. We call $\omega_k$ the \emph{rounding
term}. By Lemma~\ref{lem:backward-error}(ii),
$\tilde\eta_k = \|\widetilde A\widetilde M_k^{-1} - I\|_2$ is the
deviation~\eqref{eq:eta-def} for $\widetilde A$. How far it
is from $\eta_k$ depends on the sensitivity of the incomplete factors to
perturbations of $A$, which we do not analyse.

\begin{corollary}
\label{cor:main}
Define the growth factors
\begin{equation}
  \theta_k
  :=
  \frac{\bigl\|\,|\widetilde L_k|\,|\widetilde U_k|\,\bigr\|_2}{\|A\|_2},
  \qquad
  \Gamma(k)
  :=
  \frac{\|\widetilde L_k\|_F\,\|\widetilde U_k\|_F}{\|A\|_2},
  \label{eq:growth-def}
\end{equation}
the rounding parameters
\begin{equation}
  \mu_k^{\theta} := u\,\theta_k\,\|A\|_2\,\|\widetilde M_k^{-1}\|_2,
  \qquad
  \mu_k := u\,\Gamma(k)\,\|A\|_2\,\|\widetilde M_k^{-1}\|_2 ,
  \label{eq:stability-param}
\end{equation}
and $C_k := \gamma_{m_k}/u = m_k/(1 - m_ku)$. Under the assumptions of
Theorem~\ref{thm:stability}, with $\tilde\eta_k$ and $\omega_k$ as
in~\eqref{eq:main-bound-revised},
$\omega_k \le C_k\,\mu_k^{\theta} \le C_k\,\mu_k$ and
\[
  \|A\widetilde{M}^{-1}_k - I\|_2
  \leq
  \tilde\eta_k
  +
  C_k\,\mu_k^{\theta}
  \leq
  \tilde\eta_k
  +
  C_k\,\mu_k .
\]
\end{corollary}

\begin{proof}
Substitute~\eqref{eq:growth-def} and~\eqref{eq:stability-param}
into~\eqref{eq:rho-bound}, and use~\eqref{eq:main-bound-revised}.
\end{proof}

The growth factors and rounding parameters are those of the computed factors,
which are what is computed in practice.
$\theta_k$ comes from taking norms in Lemma~\ref{lem:backward-error}(i), and
$\Gamma(k)$ from then separating $|\widetilde L_k|$ from $|\widetilde U_k|$.
Since $\theta_k \le \Gamma(k)$, the separation can lose a great deal, but
$\|\widetilde L_k\|_F\|\widetilde U_k\|_F$ is accumulated during the
factorization at no cost, whereas $\theta_k$ needs a power method in which each
iteration applies $|\widetilde U_k|$ and then $|\widetilde L_k|$. Both growth
factors, and $\mu_k^\theta$ and $\mu_k$, are unchanged under $A \mapsto cA$
with $c > 0$, since the factors become, in exact arithmetic, $L_k$ and
$cU_k$ (for the computed factors, exactly if $c$ is a power of $2$).

\subsection{Interpretation}
\label{sec:interpretation}

Theorem~\ref{thm:stability} bounds the deviation of the preconditioned matrix
from the identity by the sum of two terms. The first, $\tilde\eta_k$, comes
from the entries that the incomplete factorization discards. The second, the
rounding term $\omega_k$, is at most $C_k\mu_k$ by Corollary~\ref{cor:main},
where $C_k \approx m_k$ and
\[
  \mu_k = u\,\|\widetilde L_k\|_F\,\|\widetilde U_k\|_F\,\|\widetilde M_k^{-1}\|_2
\]
is the unit roundoff times the size of the factors times the size of the
inverse of the preconditioner ($\mu_k^\theta \le \mu_k$ is a sharper version).
If $m_k\mu_k \ll 1$, the rounding term is small whatever the value of
$\tilde\eta_k$; otherwise the bound gives no such guarantee. In
fp64, $\mu_k$ is small unless the factors grow or $\|\widetilde M_k^{-1}\|_2$
is very large, and the bound then reduces to $\tilde\eta_k$. In lower
precision the rounding term can dominate even for moderate growth.

The bounds are upper bounds and do not predict how well the preconditioner
performs. $\mu_k \gtrsim 1$ does not imply that convergence degrades: a large
$\|\widetilde M_k^{-1}\|_2$ reflects worst-case sensitivity, whereas Krylov
convergence depends on properties of $A\widetilde M_k^{-1}$, such as
eigenvalue clustering and non-normality, that $\|\widetilde M_k^{-1}\|_2$ does
not determine.

\subsection{Illustrative examples for ILU($k$)}
\label{sec:LU examples}
The test matrices, from the SuiteSparse Collection\footnote{\url{https://sparse.tamu.edu/}},
are listed in Table~\ref{T:test problems}; they are chosen to illustrate
different behaviours. The right-hand side $b$ is chosen so that the solution is
the vector of ones. For the SPD problems, $\kappa_2(\widehat A)$ is obtained
by Lanczos on $\widehat A$ and on $\widehat A^{-1}$, the latter applied through
a sparse LU factorization; for the nonsymmetric problems, it is obtained in the
same way from $\widehat A^T\widehat A$ and its inverse.

\begin{table}[htbp]
\caption{Test problems. $n$ is the dimension of $A$ and $nnz(A)$ the number of
entries of $A$ (nonsymmetric problems) or of its lower triangle (SPD problems).
$\kappa_2(\widehat A)$ is the condition number of the scaled matrix
$\widehat A$ of Section~\ref{sec:scaling}, to which the factorization is
applied.}
\label{T:test problems}\vspace{3mm}
{
\footnotesize
\begin{center}
\begin{tabular}{lrrl} \hline
{\Hstrut Identifier} &
{$n$} &
 {$nnz(A)$} &
 {$\kappa_2(\widehat A)$} \\
 \hline
\multicolumn{4}{c}{\bf Nonsymmetric problems} \\
\hline\Tstrut
{\tt HB/sherman5}              &  3312  &   2.08$\times 10^4$ &        2.23$\times 10^3$ \\
{\tt HB/young3c}               &  841   &   3.99$\times 10^3$ &        2.41$\times 10^3$ \\
{\tt Bai/af23560}              &  23560 &   4.84$\times 10^5$ &        7.95$\times 10^3$ \\
{\tt Bomhof/circuit\_3}        &  12127 &   4.81$\times 10^4$ &        4.86$\times 10^6$ \\
\hline
\multicolumn{4}{c}{\bf SPD problems} \\
\hline\Tstrut
{\tt Nasa/nasa2146}            &  2146 &   3.72$\times 10^4$ &        9.34$\times 10^2$ \\
{\tt Boeing/msc01440}          &  1440 &   2.39$\times 10^4$ &        1.27$\times 10^{4}$ \\
{\tt HB/bcsstk24}              &  3562 &   8.17$\times 10^4$ &        1.34$\times 10^{7}$ \\
{\tt Cylshell/s2rmt3m1}        &  5489 &   1.12$\times 10^5$ &        1.77$\times 10^{7}$ \\
{\tt Cylshell/s3rmt3m1}        &  5489 &   1.12$\times 10^5$ &        1.51$\times 10^{9}$ \\
{\tt Williams/cant}            & 62451 &   2.03$\times 10^6$ &        5.59$\times 10^{10}$ \\
{\tt Boeing/msc01050}          &  1050 &   1.51$\times 10^4$ &        5.60$\times 10^{13}$ \\
\hline
\end{tabular}
\end{center}
}
\end{table}

All computation in this section is in fp64, and all quantities are those of
the scaled matrix $\widehat A$ and its computed factors $\widetilde L_k$ and
$\widetilde U_k$. We report the measured deviation
$\|\widehat A\widetilde M_k^{-1} - I\|_2$ as $\eta_k$. This assumes that
$\omega_k$ and the difference between $\tilde\eta_k$ and $\eta_k$ are
negligible, which can fail when $\widehat A$ is very ill conditioned.
$\|\widehat A\|_2$, $\|\widetilde M_k^{-1}\|_2$, $\|\widetilde R_k\|_2$ (with
$\widetilde R_kw = \widehat Aw - \widetilde L_k(\widetilde U_kw)$) and $\eta_k$ are
estimated by the power method, and $\Gamma(k)$ and $\mu_k$ are formed from
these estimates and from $\|\widetilde L_k\|_F$ and $\|\widetilde U_k\|_F$.
GMRES, with right preconditioning and without restarts, stops when the relative
residual is below $10^{-12}$, as for the conjugate gradient (CG) method in
Section~\ref{sec:experiments}, or after $800$ iterations. No pivot is
replaced in these runs.

For {\tt sherman5} (Table~\ref{T:ilu}), the iteration count decreases as $k$
increases. $\Gamma(k)$ is almost constant, and $\|\widetilde R_k\|_2$ decreases while
$\|\widetilde M_k^{-1}\|_2$ increases. $\mu_k$ is negligible, so the
bound~\eqref{eq:main-bound-revised} reduces to the deviation. $\eta_k$ stays
close to $1$, and the product $\|\widetilde R_k\|_2\,\|\widetilde M_k^{-1}\|_2$
overestimates it by a factor of about $3$ to $4$.
\begin{table}[htbp]
\caption{ILU($k$) in fp64 for the four nonsymmetric problems. All quantities are for the scaled matrix $\widehat A$ and its computed factors; $\eta_k$ is the measured $\|\widehat A\widetilde M_k^{-1} - I\|_2$. $fill$ is the number of entries in $\widetilde L_k + \widetilde U_k$, $\Gamma(k)$ and $\mu_k$ are as in \eqref{eq:growth-def} and \eqref{eq:stability-param}, $\mathrm{condest}$ is as in \eqref{eq:condest-def}, and $its$ is the number of GMRES iterations (\ddag: no convergence in $800$ iterations).}
\label{T:ilu}
{\footnotesize
\begin{center}
\tabcolsep=0.9mm
\begin{tabular}{rlllllllr} \hline
\multicolumn{1}{c}{\Hstrut $k$} & \multicolumn{1}{c}{$fill$} & \multicolumn{1}{c}{$\Gamma(k)$} & \multicolumn{1}{c}{$\|\widetilde M_k^{-1}\|_2$} & \multicolumn{1}{c}{$\mathrm{condest}$} & \multicolumn{1}{c}{$\|\widetilde R_k\|_2$} & \multicolumn{1}{c}{$\eta_k$} & \multicolumn{1}{c}{$\mu_k$} & \multicolumn{1}{c}{$its$} \\
\hline
\multicolumn{9}{l}{\Tstrut{\tt sherman5} ($n = 3312$, $\kappa_2(\widehat A) = 2.23\times 10^{3}$)} \\
0 & $2.08\times 10^{4}$ & $1.33\times 10^{3}$ & $9.08$ & $1.00$ & $3.84\times 10^{-1}$ & $1.11$ & $4.41\times 10^{-12}$ & $38$ \\
1 & $3.75\times 10^{4}$ & $1.32\times 10^{3}$ & $2.29\times 10^{1}$ & $1.64$ & $1.59\times 10^{-1}$ & $1.13$ & $1.11\times 10^{-11}$ & $25$ \\
2 & $6.39\times 10^{4}$ & $1.32\times 10^{3}$ & $3.19\times 10^{1}$ & $2.23$ & $9.89\times 10^{-2}$ & $1.14$ & $1.54\times 10^{-11}$ & $20$ \\
3 & $1.06\times 10^{5}$ & $1.32\times 10^{3}$ & $5.38\times 10^{1}$ & $3.67$ & $6.71\times 10^{-2}$ & $1.20$ & $2.60\times 10^{-11}$ & $16$ \\
5 & $2.21\times 10^{5}$ & $1.32\times 10^{3}$ & $1.36\times 10^{2}$ & $8.91$ & $3.35\times 10^{-2}$ & $1.23$ & $6.56\times 10^{-11}$ & $11$ \\
10 & $5.52\times 10^{5}$ & $1.32\times 10^{3}$ & $4.36\times 10^{2}$ & $2.98\times 10^{1}$ & $7.15\times 10^{-3}$ & $9.09\times 10^{-1}$ & $2.11\times 10^{-10}$ & $7$ \\
\hline
\multicolumn{9}{l}{\Tstrut{\tt young3c} ($n = 841$, $\kappa_2(\widehat A) = 2.41\times 10^{3}$)} \\
0 & $3.99\times 10^{3}$ & $7.68\times 10^{6}$ & $6.06\times 10^{4}$ & $3.22\times 10^{2}$ & $3.54$ & $1.24\times 10^{4}$ & $1.18\times 10^{-4}$ & $236$ \\
1 & $5.42\times 10^{3}$ & $1.20\times 10^{4}$ & $6.21\times 10^{7}$ & $1.88\times 10^{4}$ & $5.94\times 10^{1}$ & $6.97\times 10^{6}$ & $1.89\times 10^{-4}$ & $448$ \\
2 & $6.76\times 10^{3}$ & $6.15\times 10^{4}$ & $6.85\times 10^{7}$ & $1.24\times 10^{6}$ & $1.57\times 10^{2}$ & $1.54\times 10^{7}$ & $1.07\times 10^{-3}$ & $596$ \\
3 & $9.18\times 10^{3}$ & $3.21\times 10^{6}$ & $8.61\times 10^{5}$ & $5.01\times 10^{3}$ & $2.36\times 10^{1}$ & $2.29\times 10^{5}$ & $7.01\times 10^{-4}$ & $473$ \\
5 & $1.33\times 10^{4}$ & $5.44\times 10^{4}$ & $4.98\times 10^{4}$ & $6.64\times 10^{1}$ & $6.13$ & $8.56\times 10^{3}$ & $6.87\times 10^{-7}$ & $197$ \\
10 & $2.06\times 10^{4}$ & $3.69\times 10^{4}$ & $3.18\times 10^{4}$ & $8.84\times 10^{1}$ & $9.40$ & $6.02\times 10^{3}$ & $2.97\times 10^{-7}$ & $168$ \\
\hline
\multicolumn{9}{l}{\Tstrut{\tt af23560} ($n = 23560$, $\kappa_2(\widehat A) = 7.95\times 10^{3}$)} \\
0 & $4.84\times 10^{5}$ & $8.48\times 10^{3}$ & $1.14\times 10^{15}$ & $3.73\times 10^{12}$ & $1.40\times 10^{2}$ & $1.21\times 10^{15}$ & $6.47\times 10^{3}$ & \ddag \\
1 & $1.13\times 10^{6}$ & $6.24\times 10^{3}$ & $5.50\times 10^{2}$ & $4.68\times 10^{1}$ & $9.32\times 10^{-1}$ & $5.41\times 10^{1}$ & $2.31\times 10^{-9}$ & $56$ \\
2 & $2.37\times 10^{6}$ & $6.26\times 10^{3}$ & $1.11\times 10^{3}$ & $1.32\times 10^{2}$ & $5.53\times 10^{-1}$ & $1.35\times 10^{2}$ & $4.66\times 10^{-9}$ & $37$ \\
3 & $3.69\times 10^{6}$ & $6.29\times 10^{3}$ & $1.35\times 10^{3}$ & $7.52\times 10^{1}$ & $3.85\times 10^{-1}$ & $7.25\times 10^{1}$ & $5.70\times 10^{-9}$ & $23$ \\
5 & $6.11\times 10^{6}$ & $6.31\times 10^{3}$ & $1.53\times 10^{3}$ & $5.66\times 10^{1}$ & $2.62\times 10^{-1}$ & $3.77\times 10^{1}$ & $6.50\times 10^{-9}$ & $15$ \\
10 & $1.08\times 10^{7}$ & $6.34\times 10^{3}$ & $1.38\times 10^{3}$ & $5.86\times 10^{1}$ & $1.48\times 10^{-1}$ & $1.35\times 10^{1}$ & $5.88\times 10^{-9}$ & $10$ \\
\hline
\multicolumn{9}{l}{\Tstrut{\tt circuit\_3} ($n = 12127$, $\kappa_2(\widehat A) = 4.86\times 10^{6}$)} \\
0 & $4.83\times 10^{4}$ & $1.29\times 10^{11}$ & $2.01\times 10^{5}$ & $4.97\times 10^{5}$ & $1.08\times 10^{3}$ & $7.57\times 10^{4}$ & $7.01\times 10^{1}$ & \ddag \\
1 & $5.69\times 10^{4}$ & $1.28\times 10^{11}$ & $4.00\times 10^{6}$ & $7.29\times 10^{6}$ & $1.13\times 10^{3}$ & $4.04\times 10^{6}$ & $1.39\times 10^{3}$ & \ddag \\
2 & $6.17\times 10^{4}$ & $1.28\times 10^{11}$ & $2.57\times 10^{5}$ & $1.53\times 10^{6}$ & $4.18\times 10^{1}$ & $2.40\times 10^{3}$ & $8.91\times 10^{1}$ & $335$ \\
3 & $6.38\times 10^{4}$ & $1.28\times 10^{11}$ & $1.80\times 10^{6}$ & $7.77\times 10^{6}$ & $1.37$ & $1.21\times 10^{3}$ & $6.26\times 10^{2}$ & $53$ \\
5 & $6.55\times 10^{4}$ & $1.29\times 10^{11}$ & $1.77\times 10^{6}$ & $7.50\times 10^{6}$ & $6.31$ & $1.56\times 10^{3}$ & $6.16\times 10^{2}$ & $39$ \\
10 & $6.80\times 10^{4}$ & $1.28\times 10^{11}$ & $1.77\times 10^{6}$ & $5.34\times 10^{6}$ & $7.27\times 10^{1}$ & $2.98\times 10^{3}$ & $6.14\times 10^{2}$ & $15$ \\
\hline
\end{tabular}
\end{center}
}
\end{table}

For {\tt young3c}, $\Gamma(k)$ varies by orders of magnitude, and $\Gamma(k)$,
$\|\widetilde M_k^{-1}\|_2$ and the iteration count all vary non-monotonically
with $k$.

For {\tt af23560}, ILU(0) is ineffective: $\|\widetilde M_0^{-1}\|_2$, and with it
$\eta_0$ and $\mu_0$, are very large, and GMRES does not converge. For
$k \ge 1$, $\mu_k$ is small and GMRES converges, although $\Gamma(k)$ changes
little from $k = 0$. $\Gamma(k)$ alone is therefore not an indicator of
preconditioner quality.

For {\tt circuit\_3}, $\Gamma(k)$ and $\|\widetilde M_k^{-1}\|_2$ are very large, so
$\mu_k \gg 1$ at every $k$ and the bound gives no information. GMRES converges
for $k \ge 2$, quickly from $k = 3$, but not for $k \le 1$. A small $\mu_k$ is
therefore not necessary for the preconditioner to perform well. Iteration
counts in the hundreds are sensitive to rounding: with the more recent reference BLAS
distributed with LAPACK~\cite{lapack}, which computes vector norms differently, the {\tt young3c} counts for $k \ge 1$ differ by up to a third, and
for {\tt circuit\_3} at $k = 0$ GMRES converges in $406$ iterations instead of
stagnating just above the tolerance.

To detect ill-conditioned triangular factors, Chow and Saad~\cite{chsa:97}
propose monitoring
\begin{equation}
  \mathrm{condest} := \|\widetilde M_k^{-1}e\|_\infty ,
  \label{eq:condest-def}
\end{equation}
where $e$ is the
vector of all ones, which is less costly than estimating $\|\widetilde M_k^{-1}\|_2$.
Since $\|e\|_\infty = 1$, $\mathrm{condest}$ is a lower bound on
$\|\widetilde M_k^{-1}\|_\infty$ obtained from one solve with $\widetilde M_k$.
In Table~\ref{T:ilu}, $\mathrm{condest}$ follows the large changes in
$\|\widetilde M_k^{-1}\|_2$ ({\tt af23560} at $k = 0$, {\tt young3c} at $k = 1$ and $2$),
but can differ from it by more than three orders of magnitude ({\tt young3c}).

Except for {\tt sherman5} at $k = 10$, $\eta_k > 1$, but GMRES converges in
most cases: $\eta_k < 1$ is sufficient for convergence, but not necessary. None of $\eta_k$, $\Gamma(k)$ and $\mu_k$ predicts the iteration
count: they can be large while GMRES converges quickly ({\tt circuit\_3}, $k \ge 3$),
and $\Gamma(k)$ can change little while the iteration count changes greatly
({\tt af23560} from $k = 0$ to $1$).

\section{Incomplete Cholesky factorization IC($k$)}
\label{sec:ic}

In this section $A$ is SPD and the preconditioned system is solved by CG.
The IC($k$) factorization produces a lower triangular factor $L_k$ such
that
\[
  A = L_k L_k^T + R_k = M_k + R_k,
\]
where $R_k$ is the factorization residual determined
by the level-of-fill structure.
In floating-point arithmetic we compute a factor $\widetilde{L}_k$, and the
preconditioner of~\eqref{eq:M-def} becomes
$\widetilde{M}_k = \widetilde{L}_k \widetilde{L}_k^T$.

For IC($k$), $\widetilde U_k = \widetilde L_k^T$. The off-diagonal entries of
$\widetilde L_k$ are computed as in Lemma~\ref{lem:backward-error}; each
diagonal entry $\tilde\ell_{ii} = \mathrm{fl}(\tilde x_i^{1/2})$, with
$\tilde x_i = \mathrm{fl}(a_{ii} - \sum_{j<i}\tilde\ell_{ij}^2)$, has one extra
rounding, from the square root. The bound for $\tilde x_i$ (with $\gamma_{p-1}$,
as there is no division) and $\tilde\ell_{ii}^2 = \tilde x_i(1+\delta)^2$,
$|\delta| \le u$, together give Lemma~\ref{lem:backward-error} with
$\gamma_{m_k+1}$ in place of $\gamma_{m_k}$ (see also~\cite[Thm.~10.3]{high:2002}).

In exact arithmetic $\Delta A_k = 0$, so by Lemma~\ref{lem:backward-error}(ii)
the product of the incomplete factors reproduces the entries of $A$ within the
pattern (see also~\cite{meva:77,sctu:2023}). If $L$ is an incomplete Cholesky
factor of $A$, in exact arithmetic and without breakdown, under a pattern
containing the diagonal, and $R = A - LL^T$, then
\begin{equation}
  \mathrm{diag}(LL^T) = \mathrm{diag}(A),
  \qquad\text{equivalently}\qquad
  \mathrm{diag}(R) = 0.
  \label{eq:diag-exact}
\end{equation}
This fails for variants that modify the diagonal, such as modified incomplete
Cholesky with diagonal compensation. Taking traces gives the following result,
which has no ILU($k$) analogue, since $\mathrm{diag}(L_kU_k) = \mathrm{diag}(A)$
does not determine $\|L_k\|_F$ or $\|U_k\|_F$.

\begin{lemma}[Trace identity]
\label{lem:trace}
Let $A$ be SPD.
\begin{enumerate}
\item[\textup{(i)}] If $L$ is the incomplete Cholesky factor of $A$, in exact
  arithmetic and without breakdown, under any pattern containing the diagonal,
  then
  \[
    \|L\|_F^2 = \mathrm{tr}(A),
  \]
  independently of the choice of pattern.
\item[\textup{(ii)}] If the IC($k$) factor $\widetilde L_k$ is computed in
  floating-point arithmetic without breakdown, and $\gamma_{m_k+1} < 1$, then
  \begin{equation}
    \frac{\mathrm{tr}(A)}{1 + \gamma_{m_k+1}}
    \;\le\; \|\widetilde L_k\|_F^2 \;\le\;
    \frac{\mathrm{tr}(A)}{1 - \gamma_{m_k+1}} .
    \label{eq:trace-computed}
  \end{equation}
\end{enumerate}
\end{lemma}

\begin{proof}
(i) is the trace of~\eqref{eq:diag-exact}. (ii) Since the diagonal lies in
$\mathcal S$, $\mathrm{tr}(\widetilde R_k) = \mathrm{tr}(\Delta A_k)$, so
$\|\widetilde L_k\|_F^2 = \mathrm{tr}(A) - \mathrm{tr}(\Delta A_k)$. By
Lemma~\ref{lem:backward-error}(i) in its IC($k$) form,
$|\mathrm{tr}(\Delta A_k)| \le \gamma_{m_k+1}
\mathrm{tr}(|\widetilde L_k||\widetilde L_k^T|) =
\gamma_{m_k+1}\|\widetilde L_k\|_F^2$, which gives~\eqref{eq:trace-computed}.
\end{proof}

Because $\mathrm{diag}(R_k) = 0$, the entries in row $i$ of the factor satisfy
\begin{equation}
  \ell_{ii}^2 = a_{ii} - \sum_{j<i} \ell_{ij}^2 ,
  \label{eq:pivot-row}
\end{equation}
so the pivot in row $i$ is fixed by the sum of squares
$\sum_{j<i}\ell_{ij}^2$ of the off-diagonal entries in that row of $L_k$, and
if enlarging the pattern increases this sum then $\ell_{ii}$ decreases.
A single small pivot makes the preconditioner ill conditioned: $L_k^{-1}$ is
lower triangular with diagonal entries $\ell_{ii}^{-1}$, so
\begin{equation}
  \|M_k^{-1}\|_2 = \|L_k^{-1}\|_2^2 \ge \frac{1}{\min_i \ell_{ii}^2} .
  \label{eq:pivot-bound}
\end{equation}
Increasing $k$ can therefore reduce $\|R_k\|_2$ while degrading the
conditioning of $M_k$. A nearly singular $M_k$ can still be effective if its
near-zero eigenvalue corresponds to a near-zero eigenvalue of $A$.

When IC($k$) breaks down, the factors of $A + \alpha I$ ($\alpha > 0$) are used
to precondition $A$; strategies for avoiding breakdown are discussed
in~\cite{sctu:2024}. The preconditioner typically deteriorates as $\alpha$
increases.
Applied to $A + \alpha I$, \eqref{eq:diag-exact} and
Lemma~\ref{lem:trace} give, in exact arithmetic,
\begin{equation}
  \|L_k\|_F^2 = \mathrm{tr}(A + \alpha I) = \mathrm{tr}(A) + \alpha\,n .
  \label{eq:trace-shift}
\end{equation}

\subsection{Deviation bound for IC(\textit{k})}
\label{sec:ic-theorem}

\begin{corollary}[Deviation bound for IC($k$)]
\label{thm:ic-stability}
Let $A$ be SPD, and let the IC($k$) factorization of $A$ be computed in
floating-point arithmetic without breakdown, with $\gamma_{m_k+1} < 1$, giving
$\widetilde M_k = \widetilde L_k\widetilde L_k^T$. Then
Theorem~\ref{thm:stability} holds with $\widetilde U_k = \widetilde L_k^T$ and
$\gamma_{m_k+1}$ in place of $\gamma_{m_k}$.
Moreover, defining the growth factor
$\theta_k^{\mathrm{IC}} := \||\widetilde L_k|\|_2^2/\|A\|_2$ and the IC
rounding parameters
\[
  \mu_k^{\mathrm{IC},\theta}
    := u\,\theta_k^{\mathrm{IC}}\,\|A\|_2\,\|\widetilde M^{-1}_k\|_2,
  \qquad
  \mu_k^{\mathrm{IC}}
    := u\,\mathrm{tr}(A)\,\|\widetilde M^{-1}_k\|_2,
\]
and $C_k^{\mathrm{IC}} := \gamma_{m_k+1}/u \approx m_k + 1$,
\begin{equation}
  \|A\widetilde{M}^{-1}_k - I\|_2
    \leq
    \tilde\eta_k + C_k^{\mathrm{IC}}\,\mu_k^{\mathrm{IC},\theta}
    \leq
    \tilde\eta_k + \frac{C_k^{\mathrm{IC}}}{1 - \gamma_{m_k+1}}\,\mu_k^{\mathrm{IC}}.
  \label{eq:ic-main-approx}
\end{equation}
\end{corollary}

\begin{proof}
The first statement follows from the IC($k$) form of
Lemma~\ref{lem:backward-error}, part~(i) of which gives
\[
  \|\Delta A_k\|_2 \le \gamma_{m_k+1}\||\widetilde L_k||\widetilde L_k^T|\|_2
  = \gamma_{m_k+1}\||\widetilde L_k|\|_2^2 ,
\]
so
\[
  \omega_k \le \|\Delta A_k\|_2\|\widetilde M_k^{-1}\|_2 \le
  \gamma_{m_k+1}\||\widetilde L_k|\|_2^2\|\widetilde M_k^{-1}\|_2 =
  C_k^{\mathrm{IC}}\,\mu_k^{\mathrm{IC},\theta},
\]
and
by~\eqref{eq:trace-computed},
$\||\widetilde L_k|\|_2^2 \le \|\widetilde L_k\|_F^2 \le
\mathrm{tr}(A)/(1-\gamma_{m_k+1})$.
\end{proof}

When IC($k$) breaks down and $A + \alpha I$ is factorized instead,
Lemma~\ref{lem:backward-error} and Corollary~\ref{thm:ic-stability} hold with
$A + \alpha I$ in place of $A$. Since $A\widetilde M_k^{-1} - I =
\bigl((A+\alpha I)\widetilde M_k^{-1} - I\bigr) - \alpha\widetilde M_k^{-1}$,
\[
  \|A\widetilde M_k^{-1} - I\|_2 \le \tilde\eta_k + \omega_k
  + \alpha\,\|\widetilde M_k^{-1}\|_2 .
\]
If $\|\widetilde M_k^{-1}\|_2 \ge \|(A+\alpha I)^{-1}\|_2$, the last term is at
least $\alpha/(\lambda_{\min}(A)+\alpha)$, which is close to $1$ when
$\alpha \gg \lambda_{\min}(A)$; the bound then gives no information.

Since $\mathrm{tr}(A) \ge \|A\|_2$, if $\lambda_{\min}(\widetilde M_k) \le
\lambda_{\min}(A)$, that is, $\|\widetilde M^{-1}_k\|_2 \ge \|A^{-1}\|_2$, then
\[
  \mu_k^{\mathrm{IC}}
  = u\,\mathrm{tr}(A)\,\|\widetilde M^{-1}_k\|_2
  \ge u\,\|A\|_2\,\|A^{-1}\|_2
  = u\,\kappa(A).
\]
This lower bound is independent of $k$ and is large if $A$ is ill conditioned.
It limits the analysis, not the preconditioner: once
$\kappa(A) \gtrsim u^{-1}$, \eqref{eq:ic-main-approx} provides no information for any $k$,
but IC($k$) can still be effective (Section~\ref{sec:experiments}).
For M-matrices this condition fails, except possibly with equality: in exact arithmetic
$R_k \le 0$ entrywise, so $A = M_k + R_k$ is a regular splitting,
$A^{-1} \ge M_k^{-1} \ge 0$ entrywise, and
$\|M^{-1}_k\|_2 \le \|A^{-1}\|_2$~\cite{meva:77}.

The sharper $\mu_k^{\mathrm{IC},\theta}$, which avoids the Frobenius norm and
the trace identity, has a similar lower bound.
Since $\bigl\||\widetilde L_k|\bigr\|_2 \ge \|\widetilde L_k\|_2$ and
$\|\widetilde L_k\|_2^2 = \|\widetilde M_k\|_2$,
\[
  \mu_k^{\mathrm{IC},\theta} \;\ge\;
  u\,\|\widetilde M_k\|_2\,\|\widetilde M^{-1}_k\|_2
  \;=\; u\,\kappa(\widetilde M_k),
\]
which, unlike the lower bound $u\,\kappa(A)$ for $\mu_k^{\mathrm{IC}}$, holds for
every $\widetilde M_k$, so the bound is informative only if $\widetilde M_k$ is
well conditioned. The cause is the splitting
$\omega_k \le \|\Delta A_k\|_2\|\widetilde M_k^{-1}\|_2$: any growth factor
from $|\Delta A_k| \le \gamma_{m_k+1}|\widetilde L_k||\widetilde L_k^T|$ is at
least $\|\widetilde M_k\|_2$.

\subsection{Positive definiteness and range}
\label{sec:ic-positivity}

The IC($k$) preconditioner must be SPD to be used with CG. If the quantity
under every square root is positive, the computed pivots $\tilde\ell_{ii}$ are
positive, $\widetilde L_k$ is nonsingular, and
$\widetilde M_k = \widetilde L_k\widetilde L_k^T$ is SPD, whatever the
rounding errors. The only failure is breakdown, a nonpositive quantity under a
square root (B1 breakdown in~\cite{sctu:2024}, which also treats an unacceptably
small pivot as a breakdown). If $\ell_{ii}^2$ in~\eqref{eq:pivot-row} is small compared with
$a_{ii}$, rounding errors of order
$\gamma_{m_k+1}a_{ii}$ in forming it can make it nonpositive;
by~\eqref{eq:pivot-bound}, $\|M_k^{-1}\|_2 \ge 1/\ell_{ii}^2$ is then large.

In low precision the factor must also stay within the range of the working
format; in fp16 the range of finite numbers is $-65\,504$ to $65\,504$.
Overflow can occur when a column is divided by its pivot (B2 breakdown) or in an
update (B3 breakdown)~\cite{sctu:2024}. If the diagonal entries of
$\widehat A$ are at most $1$, the row identity~\eqref{eq:pivot-row} bounds every
entry of the factor. This holds after Jacobi scaling, and after the scaling by
the $2$-norms of the columns used in~\cite{sctu:2024}, for which
$\hat a_{ii} = a_{ii}/\|a_i\|_2$, with $a_i$ the $i$th column of $A$.

\begin{lemma}[Bounded entries after scaling]
\label{lem:bounded-entries}
Let $\widehat A$ be SPD with $\hat a_{ii} \le 1$ for all $i$, $\alpha \ge 0$, and let the IC($k$)
factorization of $\widehat A + \alpha I$ be computed in floating-point
arithmetic with unit roundoff $u$, assume $\gamma_{m_k+1} < 1$, and suppose
that row $i$ of $\widetilde L_k$ is computed without breakdown. Then
\[
  \sum_{j \le i} \tilde\ell_{ij}^2 \;\le\; \rho
  := \frac{(1+\alpha)(1+u)}{1 - \gamma_{m_k+1}},
  \qquad\text{and hence}\qquad
  |\tilde\ell_{ij}| \le \sqrt{\rho} .
\]
In exact arithmetic, $\sum_{j \le i} \ell_{ij}^2 = \hat a_{ii} + \alpha \le 1 + \alpha$.
\end{lemma}

\begin{proof}
The exact statement is~\eqref{eq:pivot-row} applied to $\widehat A + \alpha I$. In
floating-point arithmetic, entry $(i,i)$ of Lemma~\ref{lem:backward-error}(i)
in its IC($k$) form, applied to $\widehat A + \alpha I$ rounded to the working
precision, gives the following bound; it involves only rows $1$ to $i$, so it
holds once row $i$ is computed, whether or not later rows break down:
\[
  \Bigl|\mathrm{fl}(\hat a_{ii}+\alpha) - \sum_{j \le i} \tilde\ell_{ij}^2\Bigr|
  \le \gamma_{m_k+1} \sum_{j \le i} \tilde\ell_{ij}^2 .
\]
Since $\mathrm{fl}(\hat a_{ii}+\alpha) \le (\hat a_{ii}+\alpha)(1+u) \le
(1+\alpha)(1+u)$, rearranging gives the result.
\end{proof}

For fp16, $u = 2^{-11}$, and $\rho$ is close to $1+\alpha$ whenever
$m_k$ is much smaller than $2^{11}$. The lemma gives two checks that detect
breakdown before any overflow can occur. Write
$\sigma_{ij} := \hat a_{ij} - \sum_{s<j}\tilde\ell_{is}\tilde\ell_{js}$, so that
$\tilde\ell_{ij} = \sigma_{ij}/\tilde\ell_{jj}$ for $j < i$. While row $i$ is
computed, for each $j < i$ in turn, (i) before computing $\tilde\ell_{ij}$,
check that $|\sigma_{ij}| \le \sqrt{\rho}\,\tilde\ell_{jj}$, and (ii) after
computing it, check that $\hat a_{ii} + \alpha - \sum_{s \le j}\tilde\ell_{is}^2 > 0$.
If either check fails, row $i$ cannot be completed, by the lemma, and the
factorization breaks down at row $i$. If both checks pass, no quantity can
overflow. Check (ii) keeps the sum of squares of the entries already computed
in row $i$ below about $1+\alpha$, and since
$|\hat a_{ij}| \le (\hat a_{ii}\hat a_{jj})^{1/2} \le 1$, the Cauchy--Schwarz
inequality then gives $|\sigma_{ij}| \le 1+\rho$, approximately; check (i)
gives $|\tilde\ell_{ij}| \le \sqrt{\rho}$.

The safe tests of~\cite{sctu:2024} are more general: for any numbers $a$, $b$
and $c$ in the working format, they detect whether $a - bc$, or an intermediate
quantity in computing it, exceeds the largest number in that format, whatever
the scaling. After a scaling with $\hat a_{ii} \le 1$,
Lemma~\ref{lem:bounded-entries} shows that such a test can fail only through a
B1 breakdown, that is, a pivot that is too small or negative. This explains the observation in~\cite{sctu:2024} that, once a B1
breakdown is detected and the factorization is restarted with a shift, overflow
in the scaling of a column (B2 breakdown) does not occur, and is consistent
with the recommendation to scale before rounding to half
precision~\cite{hipz:2019}.

\subsection{Connection to CG convergence}
\label{sec:ic-cg}

CG convergence depends on the eigenvalues of $M_k^{-1}A$, whereas the bounds of
Section~\ref{sec:ic-theorem} involve norms that depend on the scaling and are
at least $u\,\kappa(\widetilde M_k)$. We therefore study these eigenvalues, in
exact arithmetic: they do not depend on the scaling, they lie on both sides of
$1$, and $\eta_k$ can greatly overestimate their distance from $1$.

\paragraph{Spectrum and scaling}
The residual $R_k = A - M_k$ is symmetric, and hence so is
\begin{equation}
  S_k := L_k^{-1} R_k L_k^{-T} .
  \label{eq:ic-Sk}
\end{equation}
Then
\begin{equation}
  L_k^{-1} A L_k^{-T} = I + S_k,
  \label{eq:ic-congruence}
\end{equation}
and since $M^{-1}_kA = L_k^{-T}(I + S_k)L_k^T$ is similar to $I + S_k$, the
eigenvalues of the preconditioned operator are $1 + \lambda_i(S_k)$ and
\[
  \kappa(I + S_k)
  = \frac{1 + \lambda_{\max}(S_k)}{1 + \lambda_{\min}(S_k)},
  \qquad 1 + \lambda_{\min}(S_k) > 0 .
\]
CG convergence is therefore governed by the spectrum of $S_k$. Its
distribution has been studied for model problems and for
M-matrices~\cite{axel:1992,axli:86a}. With $\widehat A = DAD$ for positive
diagonal $D$, $\widehat L_k = DL_k$ and $\widehat R_k = DR_kD$, so
$\widehat S_k = S_k$: the spectrum, and with it the convergence of CG in exact
arithmetic, does not depend on the prescaling. The same holds for measures
that depend only on the spectrum, such as Kaporin's
K-condition number~\cite{kapo:94}. By contrast, $\eta_k$,
$\theta_k^{\mathrm{IC}}$, $\mu_k^{\mathrm{IC}}$ and $\|M_k^{-1}\|_2$ depend on
$D$.

\medskip
\paragraph{Two-sided spectrum}
\begin{lemma}[Two-sided spectrum]
\label{lem:two-sided}
Let $A$ be SPD and let $M_k = L_kL_k^T$ be the IC($k$) preconditioner
in exact arithmetic, without breakdown under a pattern containing the diagonal, with
$R_k = A - M_k \neq 0$. Then $S_k$ given by \eqref{eq:ic-Sk} satisfies
\begin{equation}
  -1 < \lambda_{\min}(S_k) < 0 < \lambda_{\max}(S_k),
  \label{eq:two-sided}
\end{equation}
and the number of eigenvalues of $M_k^{-1}A$ lying below $1$ equals the number
of negative eigenvalues of $R_k$.
\end{lemma}

\begin{proof}
By~\eqref{eq:diag-exact}, $\mathrm{diag}(R_k) = 0$ and hence
$\mathrm{tr}(R_k) = 0$. Since $R_k$ is symmetric and nonzero it has at least one
nonzero eigenvalue, and as its eigenvalues sum to zero it must have both a
positive and a negative one.
By~\eqref{eq:ic-Sk}, $S_k$ is congruent to $R_k$, so by
Sylvester's law of inertia $S_k$ and $R_k$ have the same numbers of positive,
negative and zero eigenvalues; therefore
$\lambda_{\min}(S_k) < 0 < \lambda_{\max}(S_k)$.
Since $A$ is SPD, $I + S_k = L_k^{-1}AL_k^{-T}$ is SPD by~\eqref{eq:ic-congruence},
so $1 + \lambda_{\min}(S_k) > 0$, giving \eqref{eq:two-sided}.
Finally, the eigenvalues of $M_k^{-1}A$ are $1 + \lambda_i(S_k)$, so those below
$1$ correspond to the negative eigenvalues of $S_k$. By Sylvester's law of
inertia, $S_k$ has as many negative eigenvalues as $R_k$.
\end{proof}
Lemma~\ref{lem:two-sided} needs no assumption on $A$ beyond positive
definiteness: $\mathrm{diag}(R_k) = 0$ holds for every SPD $A$, ordering and
$k$ for which no shift is needed. A related result is due to
Manteuffel~\cite[Thm.~5.1]{mant:80}. Writing $A = \mathrm{diag}(A) - B$, so that
$-B$ is the off-diagonal part of $A$, he considers the incomplete factorization
of $\mathrm{diag}(A) - B/(1+\alpha)$ and shows that the preconditioned
eigenvalues lie on both sides of $1$ for all sufficiently large $\alpha$. His
assumption, that $\mathrm{diag}(A)^{-1}B$ is indefinite, holds whenever $A$ is
not diagonal, since $\mathrm{diag}(A)^{-1}B$ has zero trace.
Lemma~\ref{lem:two-sided} gives the same conclusion for every $k$, without a
shift.

By Lemma~\ref{lem:two-sided}, a bound on $\lambda_{\max}(S_k)$ alone does not
bound $\kappa(I + S_k)$,
and either end of the spectrum can make the preconditioner poor:
$1+\lambda_{\min}(S_k)$ close to zero, or a large $\lambda_{\max}(S_k)$. $\lambda_{\min}(S_k)$ is costly to
compute reliably (Section~\ref{sec:experiments}).

\medskip

\paragraph{The deviation $\eta_k$ versus $\|S_k\|_2$}
By~\eqref{eq:ic-Sk}, $R_kM^{-1}_k = L_kS_kL_k^{-1}$, so the two matrices have
the same eigenvalues. Since $S_k$ is symmetric, $\|S_k\|_2 =
\max_i|\lambda_i(S_k)|$, while the norm of the nonsymmetric $R_kM^{-1}_k$ is at
least its spectral radius and at most $\kappa(L_k)$ times it:
\[
  \|S_k\|_2 \le \eta_k \le \kappa(L_k)\,\|S_k\|_2 .
\]
$\eta_k$ can therefore be much larger than $\|S_k\|_2$ when $L_k$ is ill conditioned.
$\|S_k\|_2$ is unchanged by prescaling, whereas $\kappa(\widehat L_k)$ differs from $\kappa(L_k)$. For
nonsingular $B$, the Skeel condition number~\cite{skee:79} and its companion
with the factors in the opposite order are
\begin{equation}
  \mathrm{cond}(B) := \bigl\||B^{-1}|\,|B|\bigr\|_\infty,
  \qquad
  \mathrm{cond}^{*}(B) := \bigl\||B|\,|B^{-1}|\bigr\|_\infty .
  \label{eq:skeel-def}
\end{equation}
For positive diagonal $D$, $|(DB)^{-1}||DB| = |B^{-1}||B|$, so
$\mathrm{cond}(DB) = \mathrm{cond}(B)$, and likewise
$\mathrm{cond}^*(BD) = \mathrm{cond}^*(B)$~\cite[\S7.2]{high:2002}. Since
$\widehat L_k = DL_k$, $\mathrm{cond}(\widehat L_k) = \mathrm{cond}(L_k)$ and
$\mathrm{cond}^*(\widehat L_k^T) = \mathrm{cond}^*(L_k^T)$, so the hat may be
dropped. By a result of van der Sluis~\cite{vdsl:69},
\begin{equation}
  \mathrm{cond}(L_k) = \min_{D > 0} \; \kappa_\infty(DL_k),
  \label{eq:vdsluis}
\end{equation}
the minimum over positive diagonal $D$ being attained when the rows of $DL_k$
have equal $1$-norm. No diagonal scaling therefore reduces $\kappa_\infty(\widehat L_k)$
below $\mathrm{cond}(L_k)$. Since $\kappa_2(\widehat L_k) \ge \kappa_\infty(\widehat L_k)/n \ge
\mathrm{cond}(L_k)/n$, no diagonal scaling reduces the bound
$\kappa_2(\widehat L_k)$ on $\eta_k/\|S_k\|_2$ below $\mathrm{cond}(L_k)/n$;
$\mathrm{cond}(L_k)$ does not bound the ratio itself. Unlike $\lambda_{\min}(S_k)$,
$\mathrm{cond}(L_k)$ is inexpensive to estimate, in a few triangular solves by
the $1$-norm estimator~\cite{hage:84,high:88}. An upper bound needs only one:
the comparison matrix $\mathcal{M}(L_k)$, with $|(L_k)_{ii}|$ on the diagonal
and $-|(L_k)_{ij}|$ off it, satisfies $|L_k^{-1}| \le \mathcal{M}(L_k)^{-1}$
entrywise, so
$\mathrm{cond}(L_k) \le \|\mathcal{M}(L_k)^{-1}|L_k|e\|_\infty$, with equality
when the off-diagonal entries of $L_k$ are nonpositive, in particular when $A$
is an M-matrix.

Even in exact arithmetic, $\eta_k$ need not reflect the rate of convergence.
Since $\eta_k \ge \|S_k\|_2 \ge 1 - \lambda_{\min}(M_k^{-1}A)$, a single
eigenvalue of $M_k^{-1}A$ close to zero holds $\|S_k\|_2$ near $1$ whatever the
rest of the spectrum does, and $\eta_k$ at or above $1$. Increasing $k$ can then
improve the clustering, and with it the iteration count, while $\|S_k\|_2$ is
almost unchanged.

\medskip

\paragraph{Finite precision}
$\widetilde{M}_k$ replaces $M_k$, and the eigenvalues $1 + \lambda_i(S_k)$ are
perturbed. Corollary~\ref{thm:ic-stability} bounds the perturbation in norm,
with the scaling-dependent factor $\|\widetilde M_k^{-1}\|_2$. The eigenvalues
themselves, which determine CG convergence, are compared with those for the
nearby matrix $\widetilde A$ of Lemma~\ref{lem:backward-error}(ii) by a bound
that does not depend on the scaling when no shift is used.

\begin{lemma}[Rounding and the preconditioned eigenvalues]
\label{lem:eig-rounding}
Under the assumptions of Corollary~\ref{thm:ic-stability}, let
$\lambda_i := \lambda_i(\widetilde M_k^{-1}A)$ and
$\tilde\lambda_i := \lambda_i(\widetilde M_k^{-1}\widetilde A)$, $i = 1, \dots, n$,
each in decreasing order, with
$\widetilde A = A - \Delta A_k$ as in Lemma~\ref{lem:backward-error}(ii). Then,
for every $i$,
\begin{equation}
  |\lambda_i - \tilde\lambda_i|
  \le \gamma_{m_k+1}\,\bigl\||\widetilde L_k^{-1}|\,|\widetilde L_k|\bigr\|_2^2
  \le \gamma_{m_k+1}\,\mathrm{cond}(\widetilde L_k)\,
      \mathrm{cond}^*(\widetilde L_k^T),
  \label{eq:eig-rounding}
\end{equation}
and both bounds are unchanged when $\widetilde L_k$ is replaced by
$D\widetilde L_k$ for any positive diagonal $D$.
\end{lemma}

\begin{proof}
The $\lambda_i$ and $\tilde\lambda_i$ are the eigenvalues of the symmetric
matrices $\widetilde L_k^{-1}A\widetilde L_k^{-T}$ and
$\widetilde L_k^{-1}\widetilde A\widetilde L_k^{-T}$, which differ by
$F := \widetilde L_k^{-1}\Delta A_k\widetilde L_k^{-T}$. By Weyl's theorem, if
$B$ and $F$ are symmetric, the eigenvalues of $B+F$ and of $B$, each in
decreasing order, differ by at most $\|F\|_2$; that is, a symmetric
perturbation moves each eigenvalue by at most its norm. Hence
$|\lambda_i - \tilde\lambda_i| \le \|F\|_2$. By
Lemma~\ref{lem:backward-error}(i) in its IC($k$) form,
\[
  |F| \le |\widetilde L_k^{-1}||\Delta A_k||\widetilde L_k^{-T}|
  \le \gamma_{m_k+1}CC^T, \qquad C := |\widetilde L_k^{-1}||\widetilde L_k|,
\]
so $\|F\|_2 \le \gamma_{m_k+1}\|CC^T\|_2 = \gamma_{m_k+1}\|C\|_2^2$. Finally,
$\|C\|_2^2 \le \|C\|_1\|C\|_\infty$, with $\|C\|_\infty =
\mathrm{cond}(\widetilde L_k)$ and $\|C\|_1 = \|C^T\|_\infty =
\mathrm{cond}^*(\widetilde L_k^T)$ by~\eqref{eq:skeel-def}. Replacing
$\widetilde L_k$ by $D\widetilde L_k$ leaves $C$ unchanged.
\end{proof}

Unlike Corollary~\ref{thm:ic-stability}, \eqref{eq:eig-rounding} does not involve
$\|\widetilde M_k^{-1}\|_2$ or any other quantity that depends on the scaling,
so it has no lower bound of the form $u\,\kappa(\widetilde M_k)$; its constant
is the product of the Skeel condition numbers $\mathrm{cond}(\widetilde L_k)$
and $\mathrm{cond}^*(\widetilde L_k^T)$. When no shift is used, the
factor computed for $DAD$ is $D$ times the factor computed for $A$, exactly if
the diagonal entries of $D$ are powers of $2$ and otherwise up to rounding
errors, so the bound is then unchanged by the scaling. Like $\tilde\eta_k$ in
Theorem~\ref{thm:stability}, it compares with the incomplete factorization of
$\widetilde A$, not of $A$. Since $\mathrm{cond}^*(\widetilde L_k^T) \ge 1$, the
bound is at least $\gamma_{m_k+1}\,\mathrm{cond}(\widetilde L_k)$. It is small
in fp64 unless $\mathrm{cond}(\widetilde L_k)$ is very large, but in fp16, with
$\gamma_{m_k+1} > (m_k+1)\,2^{-11}$, it is below $1$ only if
$\mathrm{cond}(\widetilde L_k) < 2^{11}/(m_k+1)$. When the bound is
not small, the factorization may also break down, so that a shift is needed.

\subsection{Illustrative examples for IC($k$)}
\label{sec:experiments}

Unless stated otherwise, Jacobi scaling is used (Section~\ref{sec:scaling}). As
in Section~\ref{sec:LU examples}, $b$ is chosen so that the
solution is the vector of ones. CG stops when the relative residual is below
$10^{-12}$, or after $800$ iterations. Unlike~\cite{sctu:2024}, CG is applied
directly to $Ax=b$, not within iterative refinement. Except in Table~\ref{T:ic-scaling}, the
tables report $k = 0, 1, 2, 3, 5, 10$; for {\tt cant} and {\tt msc01050},
IC($k$) is the complete factorization from $k = 7$ and $k = 6$.
The shift $\alpha$ is chosen as in~\cite{sctu:2014a}.

The spectral norms are estimated by the power method, as in
Section~\ref{sec:LU examples}; $\theta_k^{\mathrm{IC}}$ uses the power method
on $|\widetilde L_k|$. $\lambda_{\min}$ and $\lambda_{\max}$ are computed by
Lanczos, $\lambda_{\max}$ as the largest eigenvalue of
$\widetilde L_k^{-1}\widehat A\widetilde L_k^{-T}$ and $\lambda_{\min}$ as the
reciprocal of the largest eigenvalue of
$\widetilde L_k^{T}\widehat A^{-1}\widetilde L_k$, with $\widehat A^{-1}$ applied
through a sparse LU factorization; we use this only to report the true values.
The power method gives $\|S_k\|_2$ but not a small $\lambda_{\min}$ reliably,
and the smallest Ritz value from CG can be much larger than $\lambda_{\min}$.
$\mathrm{cond}(\widetilde L_k)$, the Skeel condition
number~\eqref{eq:skeel-def}, and $\kappa_\infty(\widetilde L_k)$ are estimated
by Higham's $1$-norm estimator~\cite{hage:84,high:88}. As in
Section~\ref{sec:LU examples}, all quantities are those of $\widehat A$ and its
computed factor $\widetilde L_k$: $\eta_k$ is the measured
$\|\widehat A\widetilde M_k^{-1} - I\|_2$, and $\lambda_{\min}$,
$\lambda_{\max}$ and $\|S_k\|_2 = \max(1-\lambda_{\min}, \lambda_{\max}-1)$ are
those of $\widetilde L_k^{-1}\widehat A\widetilde L_k^{-T}$.

In the fp32 and fp16 runs, the IC($k$) factorization of $\widehat A$, rounded
to the working precision, is computed and stored in that precision; everything
else is in fp64, as described in Section~\ref{sec:notation}. The deviation
$\|\widehat A\widetilde M_k^{-1} - I\|_2$ is estimated as $\eta_k$ is, for the
factor computed in the working precision. Table~\ref{T:ic-mixed} gives the bound
$C_k^{\mathrm{IC}}\mu_k^{\mathrm{IC},\theta}$ on the rounding term
in~\eqref{eq:ic-main-approx}, with the unit roundoff of the working precision
and $\|\widetilde M_k^{-1}\|_2$ of the factor computed in that precision.

\paragraph{Full precision}
Table~\ref{T:ic-fp64} gives the fp64 results for the SPD problems. In every run $\|\widetilde L_k\|_F = \sqrt{n(1 + \alpha)}$ to at least
four digits (by~\eqref{eq:trace-shift}; not shown), whereas $\theta_k^{\mathrm{IC}}$ generally increases with $k$.

For {\tt nasa2146}, $\widehat A$ and
$\widetilde L_k$ are well conditioned, and all the quantities agree with CG: as
$k$ increases, $\lambda_{\min}$ increases towards $1$, and $\|S_k\|_2$,
$\eta_k$ and the iteration count decrease. $\eta_k$ is close to $\|S_k\|_2$.

\begin{table}[htbp]
\caption{IC($k$) in fp64 for the six SPD problems. All quantities are for the scaled matrix $\widehat A$ and its computed factor $\widetilde L_k$: $\lambda_{\min}$ and $\lambda_{\max}$ are the extreme eigenvalues of $\widetilde M_k^{-1}\widehat A$, $\|S_k\|_2 = \max(1-\lambda_{\min}, \lambda_{\max}-1)$, and $\eta_k$ is the measured $\|\widehat A\widetilde M_k^{-1} - I\|_2$. $fill$ is the number of entries in $\widetilde L_k$, $its$ the number of CG iterations (\ddag: no convergence in $800$ iterations), and $\alpha$ the shift.}
\label{T:ic-fp64}
{\footnotesize
\begin{center}
\tabcolsep=1.2mm
\begin{adjustbox}{max width=\textwidth}
\begin{tabular}{rlllllllllr} \hline
\multicolumn{1}{c}{\Hstrut $k$} & \multicolumn{1}{c}{$fill$} & \multicolumn{1}{c}{$\lambda_{\min}$} & \multicolumn{1}{c}{$\lambda_{\max}$} & \multicolumn{1}{c}{$\|S_k\|_2$} & \multicolumn{1}{c}{$\eta_k$} & \multicolumn{1}{c}{$\mathrm{cond}(\widetilde L_k)$} & \multicolumn{1}{c}{$\theta_k^{\mathrm{IC}}$} & \multicolumn{1}{c}{$\|\widetilde M_k^{-1}\|_2$} & \multicolumn{1}{c}{$\alpha$} & \multicolumn{1}{c}{$its$}  \\
\hline
\multicolumn{11}{l}{\Tstrut{\tt nasa2146} ($n = 2146$, $\kappa_2(\widehat A) = 9.34\times 10^{2}$)} \\
0 & $3.72\times 10^{4}$ & $2.40\times 10^{-1}$ & $1.25$ & $7.60\times 10^{-1}$ & $9.26\times 10^{-1}$ & $5.00\times 10^{1}$ & $1.19$ & $1.82\times 10^{2}$ & $0$ & $19$ \\
1 & $5.20\times 10^{4}$ & $4.09\times 10^{-1}$ & $1.14$ & $5.91\times 10^{-1}$ & $6.64\times 10^{-1}$ & $5.66\times 10^{1}$ & $1.22$ & $1.92\times 10^{2}$ & $0$ & $14$ \\
2 & $6.59\times 10^{4}$ & $6.37\times 10^{-1}$ & $1.10$ & $3.63\times 10^{-1}$ & $4.17\times 10^{-1}$ & $4.95\times 10^{1}$ & $1.23$ & $2.01\times 10^{2}$ & $0$ & $10$ \\
3 & $7.89\times 10^{4}$ & $7.42\times 10^{-1}$ & $1.08$ & $2.58\times 10^{-1}$ & $2.97\times 10^{-1}$ & $5.28\times 10^{1}$ & $1.23$ & $2.05\times 10^{2}$ & $0$ & $9$ \\
5 & $1.02\times 10^{5}$ & $8.92\times 10^{-1}$ & $1.05$ & $1.08\times 10^{-1}$ & $1.29\times 10^{-1}$ & $5.55\times 10^{1}$ & $1.23$ & $2.08\times 10^{2}$ & $0$ & $7$ \\
10 & $1.45\times 10^{5}$ & $9.87\times 10^{-1}$ & $1.01$ & $1.29\times 10^{-2}$ & $1.97\times 10^{-2}$ & $5.75\times 10^{1}$ & $1.23$ & $2.08\times 10^{2}$ & $0$ & $5$ \\
\hline
\multicolumn{11}{l}{\Tstrut{\tt msc01440} ($n = 1440$, $\kappa_2(\widehat A) = 1.27\times 10^{4}$)} \\
0 & $2.39\times 10^{4}$ & $6.10\times 10^{-3}$ & $9.76\times 10^{3}$ & $9.76\times 10^{3}$ & $4.71\times 10^{4}$ & $1.03\times 10^{4}$ & $2.22$ & $5.01\times 10^{5}$ & $1.25\times 10^{-4}$ & $280$ \\
1 & $4.98\times 10^{4}$ & $1.92\times 10^{-2}$ & $1.95$ & $9.81\times 10^{-1}$ & $1.34\times 10^{1}$ & $8.94\times 10^{1}$ & $2.53$ & $4.53\times 10^{2}$ & $1.60\times 10^{-2}$ & $93$ \\
2 & $6.05\times 10^{4}$ & $6.19\times 10^{-2}$ & $1.92$ & $9.38\times 10^{-1}$ & $1.70\times 10^{1}$ & $3.28\times 10^{2}$ & $2.62$ & $1.14\times 10^{3}$ & $4.00\times 10^{-3}$ & $58$ \\
3 & $7.60\times 10^{4}$ & $3.80\times 10^{-1}$ & $9.28$ & $8.28$ & $8.31\times 10^{1}$ & $2.17\times 10^{3}$ & $2.69$ & $1.65\times 10^{4}$ & $5.00\times 10^{-4}$ & $27$ \\
5 & $8.43\times 10^{4}$ & $8.42\times 10^{-1}$ & $1.16$ & $1.58\times 10^{-1}$ & $1.30$ & $5.38\times 10^{2}$ & $2.71$ & $3.56\times 10^{3}$ & $0$ & $10$ \\
10 & $9.45\times 10^{4}$ & $9.85\times 10^{-1}$ & $1.02$ & $1.55\times 10^{-2}$ & $3.08\times 10^{-2}$ & $5.39\times 10^{2}$ & $2.71$ & $3.50\times 10^{3}$ & $0$ & $5$ \\
\hline
\multicolumn{11}{l}{\Tstrut{\tt bcsstk24} ($n = 3562$, $\kappa_2(\widehat A) = 1.34\times 10^{7}$)} \\
0 & $8.17\times 10^{4}$ & $4.14\times 10^{-6}$ & $7.84$ & $6.84$ & $1.93\times 10^{1}$ & $4.72\times 10^{1}$ & $1.13$ & $9.16\times 10^{1}$ & $1.28\times 10^{-1}$ & \ddag \\
1 & $1.25\times 10^{5}$ & $3.68\times 10^{-3}$ & $1.84$ & $9.96\times 10^{-1}$ & $1.01\times 10^{2}$ & $5.60\times 10^{2}$ & $1.15$ & $1.05\times 10^{5}$ & $0$ & $96$ \\
2 & $1.73\times 10^{5}$ & $6.83\times 10^{-3}$ & $1.77$ & $9.93\times 10^{-1}$ & $8.16\times 10^{1}$ & $5.60\times 10^{2}$ & $1.16$ & $1.24\times 10^{5}$ & $0$ & $65$ \\
3 & $2.27\times 10^{5}$ & $1.02\times 10^{-2}$ & $1.85$ & $9.90\times 10^{-1}$ & $1.01\times 10^{2}$ & $8.46\times 10^{2}$ & $1.18$ & $1.54\times 10^{5}$ & $0$ & $51$ \\
5 & $3.52\times 10^{5}$ & $2.13\times 10^{-2}$ & $3.49$ & $2.49$ & $2.09\times 10^{2}$ & $2.02\times 10^{3}$ & $1.29$ & $3.28\times 10^{5}$ & $0$ & $38$ \\
10 & $7.26\times 10^{5}$ & $5.31\times 10^{-2}$ & $1.73$ & $9.47\times 10^{-1}$ & $1.29\times 10^{2}$ & $3.06\times 10^{3}$ & $1.40$ & $4.16\times 10^{5}$ & $0$ & $20$ \\
\hline
\multicolumn{11}{l}{\Tstrut{\tt s2rmt3m1} ($n = 5489$, $\kappa_2(\widehat A) = 1.77\times 10^{7}$)} \\
0 & $1.12\times 10^{5}$ & $1.26\times 10^{-6}$ & $1.57$ & $1.00$ & $2.71$ & $6.77\times 10^{1}$ & $1.38$ & $1.15\times 10^{2}$ & $0$ & $256$ \\
1 & $1.43\times 10^{5}$ & $4.00\times 10^{-6}$ & $1.26$ & $1.00$ & $4.97$ & $2.14\times 10^{2}$ & $1.62$ & $7.41\times 10^{2}$ & $0$ & $109$ \\
2 & $2.03\times 10^{5}$ & $9.49\times 10^{-6}$ & $1.72$ & $1.00$ & $5.83$ & $1.36\times 10^{2}$ & $1.73$ & $1.33\times 10^{3}$ & $0$ & $81$ \\
3 & $2.60\times 10^{5}$ & $1.95\times 10^{-5}$ & $1.63$ & $1.00$ & $4.77$ & $1.48\times 10^{2}$ & $1.79$ & $1.58\times 10^{3}$ & $0$ & $60$ \\
5 & $3.67\times 10^{5}$ & $7.56\times 10^{-5}$ & $1.33$ & $1.00$ & $3.60$ & $6.73\times 10^{2}$ & $1.86$ & $3.96\times 10^{3}$ & $0$ & $34$ \\
10 & $6.00\times 10^{5}$ & $1.26\times 10^{-2}$ & $1.34$ & $9.87\times 10^{-1}$ & $1.43\times 10^{1}$ & $1.76\times 10^{4}$ & $1.91$ & $5.96\times 10^{4}$ & $1.56\times 10^{-5}$ & $15$ \\
\hline
\multicolumn{11}{l}{\Tstrut{\tt cant} ($n = 62451$, $\kappa_2(\widehat A) = 5.59\times 10^{10}$)} \\
0 & $2.03\times 10^{6}$ & $2.26\times 10^{-10}$ & $1.04$ & $1.00$ & $1.27$ & $2.87\times 10^{1}$ & $2.29$ & $3.60$ & $5.12\times 10^{-1}$ & \ddag \\
1 & $4.50\times 10^{6}$ & $2.48\times 10^{-9}$ & $7.08$ & $6.08$ & $9.12$ & $9.14\times 10^{1}$ & $3.03$ & $8.62\times 10^{1}$ & $0$ & \ddag \\
2 & $7.20\times 10^{6}$ & $5.35\times 10^{-9}$ & $3.58$ & $2.58$ & $6.19$ & $2.10\times 10^{2}$ & $3.18$ & $1.32\times 10^{2}$ & $0$ & \ddag \\
3 & $9.95\times 10^{6}$ & $9.77\times 10^{-9}$ & $1.92$ & $1.00$ & $2.34$ & $2.84\times 10^{2}$ & $3.28$ & $1.87\times 10^{2}$ & $0$ & $692$ \\
5 & $1.46\times 10^{7}$ & $3.34\times 10^{-8}$ & $2.69$ & $1.69$ & $1.28\times 10^{1}$ & $1.52\times 10^{3}$ & $3.32$ & $1.69\times 10^{3}$ & $0$ & $411$ \\
10 & $1.68\times 10^{7}$ & $1.00$ & $1.00$ & $<10^{-4}$ & $4.61\times 10^{-6}$ & $5.93\times 10^{6}$ & $4.29$ & $7.38\times 10^{9}$ & $0$ & $1$ \\
\hline
\multicolumn{11}{l}{\Tstrut{\tt msc01050} ($n = 1050$, $\kappa_2(\widehat A) = 5.60\times 10^{13}$)} \\
0 & $1.51\times 10^{4}$ & $2.82\times 10^{-13}$ & $1.59$ & $1.00$ & $1.99$ & $1.33\times 10^{1}$ & $1.65$ & $8.42$ & $2.56\times 10^{-1}$ & $421$ \\
1 & $2.69\times 10^{4}$ & $5.84\times 10^{-8}$ & $1.82$ & $1.00$ & $5.59\times 10^{1}$ & $6.06\times 10^{3}$ & $1.82$ & $7.23\times 10^{5}$ & $0$ & $46$ \\
2 & $3.46\times 10^{4}$ & $1.13\times 10^{-11}$ & $1.22$ & $1.00$ & $5.21$ & $3.88\times 10^{2}$ & $1.81$ & $2.78\times 10^{2}$ & $8.00\times 10^{-3}$ & $46$ \\
3 & $3.74\times 10^{4}$ & $1.13\times 10^{-11}$ & $1.22$ & $1.00$ & $5.08$ & $3.84\times 10^{2}$ & $1.81$ & $2.72\times 10^{2}$ & $8.00\times 10^{-3}$ & $46$ \\
5 & $3.98\times 10^{4}$ & $9.94\times 10^{-1}$ & $1.00$ & $5.91\times 10^{-3}$ & $1.41\times 10^{-2}$ & $2.37\times 10^{7}$ & $1.82$ & $1.10\times 10^{13}$ & $0$ & $1$ \\
10 & $3.98\times 10^{4}$ & $9.94\times 10^{-1}$ & $1.00$ & $5.91\times 10^{-3}$ & $1.28\times 10^{-2}$ & $2.37\times 10^{7}$ & $1.82$ & $1.10\times 10^{13}$ & $0$ & $1$ \\
\hline
\end{tabular}
\end{adjustbox}
\end{center}
}
\end{table}

The
IC($k$) factorization of {\tt msc01440} breaks down for $k \le 3$,
and the size of the shift alone does not determine the quality of the
preconditioner. The shifts for $k = 0$ and $k = 3$ are of order $10^{-4}$, but
for $k = 0$, $\widetilde M_0$ is nearly singular and CG requires an order of magnitude more iterations than for
$k = 3$.

For {\tt bcsstk24}, a shift is needed only for $k = 0$.
For $k \ge 1$, $\eta_k$ overestimates
$\|S_k\|_2$ by a factor of about $100$. $\mathrm{cond}(\widetilde L_k)$ is larger still, so
by~\eqref{eq:vdsluis} no diagonal scaling makes
$\kappa_\infty(\widetilde L_k)$
small. $\|\widetilde M_k^{-1}\|_2$, and
with it $\mu_k^{\mathrm{IC}}$, increases with $k$ while the iteration count
decreases.

For {\tt s2rmt3m1}, the iteration count follows the
ratio $\lambda_{\max}/\lambda_{\min}$ and not $\|S_k\|_2$. From $k = 0$ to $k = 5$, the iteration
count decreases by a factor of more than seven, and
$\lambda_{\max}/\lambda_{\min}$ by a factor of more than sixty, while
$\|S_k\|_2$ stays at about $1$ and $\eta_k$ does not decrease.

The large problem {\tt cant} behaves like {\tt s2rmt3m1}: for $k \le 5$, $\lambda_{\min}$
is below $10^{-7}$ and CG is slow or does not converge, and only the complete factorization,   IC(7), gives rapid convergence.

The most ill-conditioned problem is {\tt msc01050}. For $k \ge 5$, IC($k$) is
complete or, at $k = 5$, lacks only $19$ of its entries, and CG converges in one
iteration. $\widetilde M_k^{-1}\widehat A$ is then the identity apart from
rounding errors of order $u\,\kappa_2(\widehat A)$, the measured $\eta_k$ and
$1-\lambda_{\min}$ are of this size, and the bound on the rounding term
exceeds $1$. At $k = 2$ and $3$, a shift is needed and $\lambda_{\min}$ equals
$\lambda_{\min}(\widehat A)/(\lambda_{\min}(\widehat A)+\alpha)$ to the
accuracy of the computation. This is much smaller than $\lambda_{\min}$ at
$k = 1$, where no shift is needed, but the iteration count is the same. CG does
not detect the small $\lambda_{\min}$: its smallest Ritz value at $k = 2$ is $4.2\times10^{-4}$.

\paragraph{Scaling}
Table~\ref{T:ic-scaling} reports, for {\tt bcsstk24}, the quantities in the bounds
with four diagonal scalings: none, Jacobi, the scaling of Lin and
Mor\'e~\cite{limo:99} (which is based on the $2$-norms of the columns), and
the equilibration scaling of {\tt MC77}~\cite{Ruiz:2001}. For $k \ge 1$ no
shift is needed, so the preconditioned operator is the same for all four
scalings (Section~\ref{sec:scaling}), but the quantities in
Table~\ref{T:ic-scaling} depend on the scaling. Of the four, Jacobi scaling gives the
smallest $\kappa_\infty(\widetilde L_k)$, within a factor of about $5$ of
$\mathrm{cond}(\widetilde L_k)$ (Table~\ref{T:ic-fp64}), its minimum over all
diagonal scalings by~\eqref{eq:vdsluis}. $\|\widetilde M_k^{-1}\|_2$ changes by
a factor of about $2\times10^9$, mainly because the unscaled entries are large,
whereas $\theta_k^{\mathrm{IC}}$ changes by a factor of at most $2.5$. Jacobi
scaling is canonical, since it gives the same $\widehat A$ whatever diagonal
scaling is applied to $A$ first, but the quantities still do not predict the
iteration count.

When a shift is needed, the scaling changes the preconditioner
(Section~\ref{sec:scaling}). Table~\ref{T:ic-s3rmt3m1} gives the CG iteration
counts for {\tt s3rmt3m1} with the four scalings. At $k = 1$ and $2$ no shift is
needed; in exact arithmetic the counts would be the same, so the differences
result from rounding errors, to which CG is sensitive on this very
ill-conditioned problem. At $k = 5$ and $10$ a shift is needed: without scaling
it is $250$ to $500$ times larger than with the other three scalings, and CG
needs between two and six times as many iterations. Jacobi scaling gives the
fewest iterations.

\begin{table}[htbp]
\caption{IC($k$) for {\tt bcsstk24} with four scalings: none ($D = I$), Jacobi, Lin--Mor\'e~\cite{limo:99} and {\tt MC77}~\cite{Ruiz:2001}. The quantities are for $\widehat A=DAD$ and its computed factor $\widetilde L_k$. The iteration counts, $\lambda_{\min}$, $\lambda_{\max}$, $\|S_k\|_2$ and $\mathrm{cond}(\widetilde L_k)$ are the same for the four scalings to at least two significant digits, and are given in Table~\ref{T:ic-fp64}.}
\label{T:ic-scaling}
{\footnotesize
\begin{center}
\renewcommand{\arraystretch}{1.3}
\begin{tabular}{lllll} \hline
\Tstrut & none & Jacobi & Lin--Mor\'e & {\tt MC77} \\
\hline
\multicolumn{5}{l}{\Tstrut $k = 1$} \\
$\eta_k$ & $1.54\times 10^{3}$ & $1.01\times 10^{2}$ & $2.62\times 10^{2}$ & $1.59\times 10^{2}$ \\
$\kappa_\infty(\widetilde L_k)$ & $3.58\times 10^{5}$ & $1.70\times 10^{3}$ & $6.75\times 10^{3}$ & $3.24\times 10^{3}$ \\
$\theta_k^{\mathrm{IC}}$ & $1.32$ & $1.15$ & $1.90$ & $2.26$ \\
$\|\widetilde M_k^{-1}\|_2$ & $8.18\times 10^{-4}$ & $1.05\times 10^{5}$ & $1.64\times 10^{6}$ & $8.94\times 10^{5}$ \\
$\mu_k^{\mathrm{IC}}$ & $1.21\times 10^{-4}$ & $4.16\times 10^{-8}$ & $3.95\times 10^{-7}$ & $1.26\times 10^{-7}$ \\
\hline
\multicolumn{5}{l}{\Tstrut $k = 10$} \\
$\eta_k$ & $8.98\times 10^{2}$ & $1.29\times 10^{2}$ & $3.28\times 10^{2}$ & $2.15\times 10^{2}$ \\
$\kappa_\infty(\widetilde L_k)$ & $3.10\times 10^{6}$ & $1.61\times 10^{4}$ & $5.94\times 10^{4}$ & $2.76\times 10^{4}$ \\
$\theta_k^{\mathrm{IC}}$ & $1.32$ & $1.40$ & $2.83$ & $3.28$ \\
$\|\widetilde M_k^{-1}\|_2$ & $2.22\times 10^{-3}$ & $4.16\times 10^{5}$ & $5.52\times 10^{6}$ & $3.76\times 10^{6}$ \\
$\mu_k^{\mathrm{IC}}$ & $3.29\times 10^{-4}$ & $1.64\times 10^{-7}$ & $1.33\times 10^{-6}$ & $5.30\times 10^{-7}$ \\
\hline
\end{tabular}
\end{center}
}
\end{table}

\begin{table}[htbp]
\caption{IC($k$) for {\tt s3rmt3m1} with four scalings: CG iteration count and, in brackets, the shift $\alpha$ where one is needed.}
\label{T:ic-s3rmt3m1}
{\footnotesize
\begin{center}
\begin{tabular}{rllll} \hline
\Tstrut $k$ & none & Jacobi & Lin--Mor\'e & {\tt MC77} \\
\hline\Tstrut
1 & 204 & 224 & 209 & 197 \\
2 & 67 & 67 & 67 & 70 \\
5 & 234 ($8.00\times 10^{-3}$) & 39 ($1.56\times 10^{-5}$) & 48 ($1.56\times 10^{-5}$) & 74 ($1.56\times 10^{-5}$) \\
10 & 160 ($4.00\times 10^{-3}$) & 33 ($1.56\times 10^{-5}$) & 45 ($1.56\times 10^{-5}$) & 67 ($1.56\times 10^{-5}$) \\
\hline
\end{tabular}
\end{center}
}
\end{table}

\paragraph{Low precision}
Table~\ref{T:ic-mixed} compares the three precisions on the six problems of
Table~\ref{T:ic-fp64}, with
$u\,\kappa_2(\widehat A)$ from $6\times10^{-5}$ to $3\times10^{6}$ in fp32 and
from $0.5$ to $3\times10^{10}$ in fp16.

In fp32, the iteration counts are close to those in fp64, even where the
bound on the rounding term exceeds $1$ and the deviation bound provides no
information, as on {\tt msc01440}, {\tt bcsstk24}, {\tt s2rmt3m1} and {\tt msc01050}. The exceptions are
the (nearly) complete factorizations of {\tt cant} and {\tt msc01050}, where CG needs more
than one iteration.

In fp16, IC($k$) breaks down on {\tt bcsstk24}, {\tt s2rmt3m1} and {\tt msc01050} at levels at
which fp64 needs no shift. The shift is then greater than $2^{-11}$, since a
smaller shift has no effect on the unit diagonal, and the shifted preconditioner
improves little with $k$ (see, for example, {\tt bcsstk24} for $k \ge 1$ and
{\tt s2rmt3m1} for $k \ge 5$). The iteration count can also increase without a
change of shift: for {\tt msc01050} at $k = 2$ and $3$, for {\tt msc01440} at
$k = 5$ and $10$, and at the complete factorization of {\tt cant}, where fp16
needs hundreds of iterations instead of one.
The bound on the rounding term is greater than $1$ in every run, and does not
distinguish the cases in which fp16 increases the iteration count from those in
which it does not. For {\tt bcsstk24} at $k = 10$, $\gamma_{m_k+1} \ge 1$, so
Corollary~\ref{thm:ic-stability} and Lemma~\ref{lem:bounded-entries} do not
apply.

\begin{table}[htbp]
\caption{IC($k$) in three precisions: CG iteration count (\ddag: no convergence in $800$ iterations), the shift $\alpha$ in fp64 and fp16, and the bound $C_k^{\mathrm{IC}}\mu_k^{\mathrm{IC},\theta}$ on the rounding term for fp32 and fp16; a dash indicates that $\gamma_{m_k+1} \ge 1$, so that Corollary~\ref{thm:ic-stability} does not apply.}
\label{T:ic-mixed}
{\scriptsize
\begin{center}
\tabcolsep=1.3mm
\begin{tabular}{lrrrrllll} \hline
\Tstrut problem & $k$ & \multicolumn{3}{c}{$its$} & \multicolumn{2}{c}{$\alpha$} & \multicolumn{2}{c}{$C_k^{\mathrm{IC}}\mu_k^{\mathrm{IC},\theta}$} \\
 & & fp64 & fp32 & fp16 & fp64 & fp16 & fp32 & fp16 \\
\hline\Tstrut
\Tstrut {\tt nasa2146} & 0 & $19$ & $19$ & $19$ & $0$ & $0$ & $1.22\times 10^{-3}$ & $1.03\times 10^{1}$ \\
 & 1 & $14$ & $14$ & $14$ & $0$ & $0$ & $1.82\times 10^{-3}$ & $1.52\times 10^{1}$ \\
 & 2 & $10$ & $10$ & $10$ & $0$ & $0$ & $2.44\times 10^{-3}$ & $2.07\times 10^{1}$ \\
 & 3 & $9$ & $9$ & $9$ & $0$ & $0$ & $3.04\times 10^{-3}$ & $2.54\times 10^{1}$ \\
 & 5 & $7$ & $7$ & $9$ & $0$ & $0$ & $4.18\times 10^{-3}$ & $3.51\times 10^{1}$ \\
 & 10 & $5$ & $5$ & $8$ & $0$ & $0$ & $5.70\times 10^{-3}$ & $4.82\times 10^{1}$ \\
\hline
\Tstrut {\tt msc01440} & 0 & $280$ & $291$ & $237$ & $1.25\times 10^{-4}$ & $5.00\times 10^{-4}$ & $8.73$ & $1.18\times 10^{2}$ \\
 & 1 & $93$ & $93$ & $94$ & $1.60\times 10^{-2}$ & $1.60\times 10^{-2}$ & $2.86\times 10^{-2}$ & $3.99\times 10^{2}$ \\
 & 2 & $58$ & $58$ & $59$ & $4.00\times 10^{-3}$ & $4.00\times 10^{-3}$ & $9.12\times 10^{-2}$ & $8.38\times 10^{2}$ \\
 & 3 & $27$ & $27$ & $29$ & $5.00\times 10^{-4}$ & $5.00\times 10^{-4}$ & $2.16$ & $6.13\times 10^{3}$ \\
 & 5 & $10$ & $10$ & $20$ & $0$ & $0$ & $5.32\times 10^{-1}$ & $7.27\times 10^{3}$ \\
 & 10 & $5$ & $5$ & $19$ & $0$ & $0$ & $7.48\times 10^{-1}$ & $1.08\times 10^{4}$ \\
\hline
\Tstrut {\tt bcsstk24} & 0 & \ddag & \ddag & \ddag & $1.28\times 10^{-1}$ & $1.28\times 10^{-1}$ & $1.89\times 10^{-3}$ & $1.68\times 10^{1}$ \\
 & 1 & $96$ & $96$ & $206$ & $0$ & $5.00\times 10^{-4}$ & $4.35$ & $2.28\times 10^{3}$ \\
 & 2 & $65$ & $65$ & $197$ & $0$ & $5.00\times 10^{-4}$ & $1.02\times 10^{1}$ & $5.44\times 10^{3}$ \\
 & 3 & $51$ & $51$ & $211$ & $0$ & $5.00\times 10^{-4}$ & $2.17\times 10^{1}$ & $1.04\times 10^{4}$ \\
 & 5 & $38$ & $37$ & $183$ & $0$ & $5.00\times 10^{-4}$ & $9.50\times 10^{1}$ & $2.46\times 10^{4}$ \\
 & 10 & $20$ & $20$ & $171$ & $0$ & $5.00\times 10^{-4}$ & $3.44\times 10^{2}$ & -- \\
\hline
\Tstrut {\tt s2rmt3m1} & 0 & $256$ & $256$ & $255$ & $0$ & $0$ & $1.13\times 10^{-3}$ & $8.95$ \\
 & 1 & $109$ & $109$ & $109$ & $0$ & $0$ & $1.06\times 10^{-2}$ & $9.21\times 10^{1}$ \\
 & 2 & $81$ & $81$ & $86$ & $0$ & $0$ & $2.82\times 10^{-2}$ & $3.36\times 10^{2}$ \\
 & 3 & $60$ & $60$ & $71$ & $0$ & $0$ & $4.45\times 10^{-2}$ & $1.37\times 10^{3}$ \\
 & 5 & $34$ & $34$ & $57$ & $0$ & $5.00\times 10^{-4}$ & $1.66\times 10^{-1}$ & $5.26\times 10^{2}$ \\
 & 10 & $15$ & $15$ & $53$ & $1.56\times 10^{-5}$ & $5.00\times 10^{-4}$ & $4.49$ & $1.70\times 10^{3}$ \\
\hline
\Tstrut {\tt cant} & 0 & \ddag & \ddag & \ddag & $5.12\times 10^{-1}$ & $5.12\times 10^{-1}$ & $1.52\times 10^{-4}$ & $1.27$ \\
 & 1 & \ddag & \ddag & \ddag & $0$ & $0$ & $1.14\times 10^{-2}$ & $8.96\times 10^{1}$ \\
 & 2 & \ddag & \ddag & \ddag & $0$ & $0$ & $3.32\times 10^{-2}$ & $2.56\times 10^{2}$ \\
 & 3 & $692$ & $692$ & $714$ & $0$ & $0$ & $7.66\times 10^{-2}$ & $4.76\times 10^{2}$ \\
 & 5 & $411$ & $409$ & $509$ & $0$ & $0$ & $7.02\times 10^{-1}$ & $8.88\times 10^{2}$ \\
 & 10 & $1$ & $16$ & $446$ & $0$ & $0$ & $4.86\times 10^{3}$ & $8.83\times 10^{2}$ \\
\hline
\Tstrut {\tt msc01050} & 0 & $421$ & $420$ & $421$ & $2.56\times 10^{-1}$ & $2.56\times 10^{-1}$ & $5.31\times 10^{-4}$ & $4.64$ \\
 & 1 & $46$ & $46$ & $53$ & $0$ & $5.00\times 10^{-4}$ & $9.80\times 10^{1}$ & $1.90\times 10^{3}$ \\
 & 2 & $46$ & $46$ & $68$ & $8.00\times 10^{-3}$ & $8.00\times 10^{-3}$ & $6.41\times 10^{-2}$ & $6.72\times 10^{2}$ \\
 & 3 & $46$ & $46$ & $107$ & $8.00\times 10^{-3}$ & $8.00\times 10^{-3}$ & $6.51\times 10^{-2}$ & $6.92\times 10^{2}$ \\
 & 5 & $1$ & $7$ & $56$ & $0$ & $5.00\times 10^{-4}$ & $2.28\times 10^{2}$ & $2.83\times 10^{3}$ \\
 & 10 & $1$ & $7$ & $56$ & $0$ & $5.00\times 10^{-4}$ & $2.28\times 10^{2}$ & $2.83\times 10^{3}$ \\
\hline
\end{tabular}
\end{center}
}
\end{table}

Table~\ref{T:ic-fp16} compares the fp16 deviation and $\lambda_{\min}$ with
fp64. The deviation does not reflect the increase in iterations caused by the
shift: on {\tt bcsstk24} for $k \ge 1$, and on {\tt s2rmt3m1} at $k = 5$ and $10$, it is
smaller than $\eta_k$, although fp16 needs more iterations; on {\tt bcsstk24} it also
decreases with $k$, while the fp16 iteration count changes little.
However,
$\lambda_{\min}$ does reflect it: at these levels it is below its fp64 value.
On {\tt bcsstk24}, for $k \ge 1$, it stays between $5.9\times10^{-4}$ and
$6.6\times10^{-4}$ in fp16, whereas in fp64 it increases from
$3.7\times10^{-3}$ to $5.3\times10^{-2}$. This is what the shift implies: even
the complete factorization of $\widehat A+\alpha I$ gives
\[
  \lambda_{\min}(\widehat M^{-1}\widehat A)
  = \frac{\lambda_{\min}(\widehat A)}{\lambda_{\min}(\widehat A)+\alpha},
\]
which does not depend on $k$ and is small when $\lambda_{\min}(\widehat A) \ll \alpha$. For {\tt msc01440} at $k = 5$ and $10$, and for {\tt cant} at
$k = 10$, fp16 needs more iterations than fp64 although no shift is used; there
the deviation is larger than $\eta_k$ and $\lambda_{\min}$ is smaller than in
fp64. For {\tt msc01050} at $k = 2$ and $3$, the deviation and $\lambda_{\min}$ are almost
unchanged from fp64 and do not reflect the increase.

\begin{table}[htbp]
\caption{IC($k$) in fp64 and fp16: CG iteration count (\ddag: no convergence in $800$ iterations), $\eta_k$ in fp64 and the measured deviation $\|\widehat A\widetilde M_k^{-1} - I\|_2$ in fp16, and $\lambda_{\min}$.}
\label{T:ic-fp16}
{\scriptsize
\begin{center}
\tabcolsep=1.3mm
\begin{tabular}{lrrrllll} \hline
\Tstrut problem & $k$ & \multicolumn{2}{c}{$its$} & $\eta_k$ & deviation & \multicolumn{2}{c}{$\lambda_{\min}$} \\
 & & fp64 & fp16 & fp64 & fp16 & fp64 & fp16 \\
\hline\Tstrut
\Tstrut {\tt nasa2146} & 0 & $19$ & $19$ & $9.26\times 10^{-1}$ & $9.10\times 10^{-1}$ & $2.40\times 10^{-1}$ & $2.39\times 10^{-1}$ \\
 & 1 & $14$ & $14$ & $6.64\times 10^{-1}$ & $6.68\times 10^{-1}$ & $4.09\times 10^{-1}$ & $4.03\times 10^{-1}$ \\
 & 2 & $10$ & $10$ & $4.17\times 10^{-1}$ & $4.35\times 10^{-1}$ & $6.37\times 10^{-1}$ & $6.24\times 10^{-1}$ \\
 & 3 & $9$ & $9$ & $2.97\times 10^{-1}$ & $3.27\times 10^{-1}$ & $7.42\times 10^{-1}$ & $7.21\times 10^{-1}$ \\
 & 5 & $7$ & $9$ & $1.29\times 10^{-1}$ & $1.99\times 10^{-1}$ & $8.92\times 10^{-1}$ & $8.51\times 10^{-1}$ \\
 & 10 & $5$ & $8$ & $1.97\times 10^{-2}$ & $1.56\times 10^{-1}$ & $9.87\times 10^{-1}$ & $9.10\times 10^{-1}$ \\
\hline
\Tstrut {\tt msc01440} & 0 & $280$ & $237$ & $4.71\times 10^{4}$ & $6.23\times 10^{1}$ & $6.10\times 10^{-3}$ & $6.02\times 10^{-3}$ \\
 & 1 & $93$ & $94$ & $1.34\times 10^{1}$ & $2.16\times 10^{1}$ & $1.92\times 10^{-2}$ & $2.01\times 10^{-2}$ \\
 & 2 & $58$ & $59$ & $1.70\times 10^{1}$ & $1.73\times 10^{1}$ & $6.19\times 10^{-2}$ & $6.21\times 10^{-2}$ \\
 & 3 & $27$ & $29$ & $8.31\times 10^{1}$ & $2.50\times 10^{1}$ & $3.80\times 10^{-1}$ & $2.43\times 10^{-1}$ \\
 & 5 & $10$ & $20$ & $1.30$ & $5.40$ & $8.42\times 10^{-1}$ & $3.81\times 10^{-1}$ \\
 & 10 & $5$ & $19$ & $3.08\times 10^{-2}$ & $5.31$ & $9.85\times 10^{-1}$ & $3.81\times 10^{-1}$ \\
\hline
\Tstrut {\tt bcsstk24} & 0 & \ddag & \ddag & $1.93\times 10^{1}$ & $2.04\times 10^{1}$ & $4.14\times 10^{-6}$ & $4.13\times 10^{-6}$ \\
 & 1 & $96$ & $206$ & $1.01\times 10^{2}$ & $3.27\times 10^{1}$ & $3.68\times 10^{-3}$ & $5.87\times 10^{-4}$ \\
 & 2 & $65$ & $197$ & $8.16\times 10^{1}$ & $1.60\times 10^{1}$ & $6.83\times 10^{-3}$ & $6.34\times 10^{-4}$ \\
 & 3 & $51$ & $211$ & $1.01\times 10^{2}$ & $1.22\times 10^{1}$ & $1.02\times 10^{-2}$ & $6.49\times 10^{-4}$ \\
 & 5 & $38$ & $183$ & $2.09\times 10^{2}$ & $8.46$ & $2.13\times 10^{-2}$ & $6.55\times 10^{-4}$ \\
 & 10 & $20$ & $171$ & $1.29\times 10^{2}$ & $7.67$ & $5.31\times 10^{-2}$ & $6.57\times 10^{-4}$ \\
\hline
\Tstrut {\tt s2rmt3m1} & 0 & $256$ & $255$ & $2.71$ & $2.63$ & $1.26\times 10^{-6}$ & $1.26\times 10^{-6}$ \\
 & 1 & $109$ & $109$ & $4.97$ & $5.58$ & $4.00\times 10^{-6}$ & $4.01\times 10^{-6}$ \\
 & 2 & $81$ & $86$ & $5.83$ & $8.65$ & $9.49\times 10^{-6}$ & $9.40\times 10^{-6}$ \\
 & 3 & $60$ & $71$ & $4.77$ & $1.85\times 10^{1}$ & $1.95\times 10^{-5}$ & $1.81\times 10^{-5}$ \\
 & 5 & $34$ & $57$ & $3.60$ & $2.34$ & $7.56\times 10^{-5}$ & $3.90\times 10^{-5}$ \\
 & 10 & $15$ & $53$ & $1.43\times 10^{1}$ & $3.12$ & $1.26\times 10^{-2}$ & $8.08\times 10^{-5}$ \\
\hline
\Tstrut {\tt cant} & 0 & \ddag & \ddag & $1.27$ & $1.27$ & $2.26\times 10^{-10}$ & $2.25\times 10^{-10}$ \\
 & 1 & \ddag & \ddag & $9.12$ & $8.28$ & $2.48\times 10^{-9}$ & $2.41\times 10^{-9}$ \\
 & 2 & \ddag & \ddag & $6.19$ & $5.19$ & $5.35\times 10^{-9}$ & $4.83\times 10^{-9}$ \\
 & 3 & $692$ & $714$ & $2.34$ & $1.63$ & $9.77\times 10^{-9}$ & $7.75\times 10^{-9}$ \\
 & 5 & $411$ & $509$ & $1.28\times 10^{1}$ & $1.68$ & $3.34\times 10^{-8}$ & $1.54\times 10^{-8}$ \\
 & 10 & $1$ & $446$ & $4.61\times 10^{-6}$ & $1.09$ & $1.00$ & $2.21\times 10^{-8}$ \\
\hline
\Tstrut {\tt msc01050} & 0 & $421$ & $421$ & $1.99$ & $2.00$ & $2.82\times 10^{-13}$ & $2.82\times 10^{-13}$ \\
 & 1 & $46$ & $53$ & $5.59\times 10^{1}$ & $3.22$ & $5.84\times 10^{-8}$ & $9.96\times 10^{-11}$ \\
 & 2 & $46$ & $68$ & $5.21$ & $5.44$ & $1.13\times 10^{-11}$ & $1.16\times 10^{-11}$ \\
 & 3 & $46$ & $107$ & $5.08$ & $5.30$ & $1.13\times 10^{-11}$ & $1.16\times 10^{-11}$ \\
 & 5 & $1$ & $56$ & $1.41\times 10^{-2}$ & $2.91$ & $9.94\times 10^{-1}$ & $9.95\times 10^{-11}$ \\
 & 10 & $1$ & $56$ & $1.28\times 10^{-2}$ & $2.91$ & $9.94\times 10^{-1}$ & $9.95\times 10^{-11}$ \\
\hline
\end{tabular}
\end{center}
}
\end{table}

\FloatBarrier
\section{Conclusions}
\label{sec:conclusions}

We have analysed the effect of rounding errors and diagonal scaling on
level-based ILU($k$) and IC($k$). Because the sparsity pattern is fixed in
advance, the deviation of the preconditioned operator from the identity is
bounded by an incompleteness term plus a rounding term from the backward error
on the pattern. When no shift or pivot modification is needed, prescaling leaves
the iterates unchanged in exact arithmetic but can change the quantities in this
bound by orders of magnitude (Table~\ref{T:ic-scaling}). The scaling changes the
preconditioner only when a shift or pivot modification is needed, and no single
scaling is the best for every problem. For SPD matrices, Jacobi scaling is
canonical, gives the unit diagonal on which the bound on the entries and the
effect of the shift depend, and gave the fewest iterations when a shift was
needed (Table~\ref{T:ic-s3rmt3m1}). The preconditioned eigenvalues and the
Skeel condition numbers avoid the choice, since they do not depend on the
scaling when no shift is used.

\emph{Nonsymmetric problems (ILU($k$)).} None of $\eta_k$, $\Gamma(k)$ and
$\mu_k$ predicts the GMRES iteration count: GMRES converges with
$\mu_k \gg 1$ and with $\eta_k > 1$, and $\Gamma(k)$ can change little while the
iteration count changes greatly. Relating scale-invariant quantities, such as
the constant in Lemma~\ref{lem:eig-rounding}, to the convergence of GMRES is a
direction for future work.

\emph{SPD problems (IC($k$)).} The normwise bound on the rounding term is at
least $u\,\kappa(\widetilde M_k)$, so it gives no guarantee when the
preconditioner is ill conditioned relative to the working precision. The
preconditioned eigenvalues, which lie on both sides of $1$ and determine CG
convergence, move under rounding, relative to those of the incomplete
factorization of a nearby matrix, by at most
$\gamma_{m_k+1}\,\mathrm{cond}(\widetilde L_k)\,\mathrm{cond}^*(\widetilde L_k^T)$
(Lemma~\ref{lem:eig-rounding}). Unlike the normwise bound, this has no lower
limit of the form $u\,\kappa(\widetilde M_k)$, and when no shift is used it
does not depend on the scaling; both Skeel condition numbers are inexpensive to
estimate. After a diagonal scaling that makes the diagonal entries at most $1$, such as
Jacobi scaling, the entries of the factor of $\widehat A + \alpha I$ are bounded
by about $\sqrt{1+\alpha}$, where $\alpha \ge 0$ is the shift ($\alpha = 0$
if none is needed). Provided $\gamma_{m_k+1} < 1$ (see Lemma~\ref{lem:bounded-entries}), the
factor used as the preconditioner therefore cannot overflow, even in fp16, and
in a factorization that breaks down, testing against this bound before each
division signals the breakdown before any overflow can occur.

In the IC($k$) experiments, the CG iteration count generally follows
$\lambda_{\min}$ and $\lambda_{\max}/\lambda_{\min}$ rather than $\eta_k$ or
$\|\widetilde M_k^{-1}\|_2$. Computing and storing the factor in fp32 changes the iteration count little except
near the complete factorization. In fp16, IC($k$) breaks down at levels at which
fp64 needs no shift; the shift, which must exceed $2^{-11}$, results in
the preconditioner improving little with $k$. This limit may be avoided by forming the
shifted diagonal in higher precision. The iteration count can also increase
without a shift, through rounding errors in the factor.

\section*{Conflict of interest}
The authors declare no potential conflict of interests.

\def\cprime{$'$} \def\cprime{$'$} \def\cprime{$'$}

\end{document}